\documentclass[11pt]{article}

\usepackage{lmodern}
\usepackage{cite}
\usepackage[T1]{fontenc}
\usepackage{amsthm}
\usepackage{amsfonts}
\usepackage{amssymb}
\usepackage{amsmath}
\usepackage{graphicx}
\usepackage{mathrsfs} % for \mathscr{F}
\usepackage{mathtools}
\usepackage{hyperref}
\usepackage{array,multirow}
\usepackage{caption}   % for tables
\usepackage{float}   % for tables
\usepackage{microtype}
\usepackage{subfiles}    %for splitting input file
\usepackage[title]{appendix}   % For nicer Appendices Headers
\usepackage{longtable}
\usepackage{needspace}
\usepackage{calc} 
\usepackage{xfrac} 
\usepackage{enumitem}

   \title{The Order of the Monster Finite Simple Group}

\author{Gerald H\"ohn\footnote{Department of Mathematics, Kansas State University, USA, {\tt gerald@monstrous-moonshine.de}}\phantom{r}
and  Martin Seysen\footnote{München, Germany, {\tt m.seysen@gmx.de}}}

\date{August 1, 2025}

   \newcommand{\Skip}[1]{}

   \newcommand{\setZ}{\mathbb{Z}}

   \newcommand {\SL}{\,{\mathop{\rm SL}}}

   \newcommand {\Aut}{{\mathop{\rm Aut}}}

   \newcommand {\ax}{\,{\mathop{\rm ax}}}
   \newcommand {\tr}{\mathop{\rm tr}}
   \newcommand {\ad}{\mathop{\rm ad}}

   \newcommand {\fatline}[1]{\noindent {\bf #1}}
   \newcommand {\proofend}{\noindent $\Box$ }

   \newtheorem {Theorem}{Theorem}[section]
   \newtheorem {Definition}[Theorem]{Definition}
   \newtheorem {Corollary}[Theorem]{Corollary}
   \newtheorem {Proposition}[Theorem]{Proposition}
   \newtheorem {Lemma}[Theorem]{Lemma}

   \numberwithin{equation}{Theorem}
 
\newcommand{\M}{\mathbb{M}}
\newcommand{\Griess}{\mathcal{B}}
\newcommand{\Golay}{{G_{24}}}
\newcommand{\F}{\mathbb{F}}
\newcommand{\Q}{\mathbb{Q}}
\newcommand{\R}{\mathbb{R}}
\newcommand{\Z}{\mathbb{Z}}

\newcommand{\Mv}{\mathbb{M}_{v^+}}

\def\voa{vertex operator algebra\ }
\def\V{V^{\natural}}
\def\Co{{\rm Co}}
\def\C{{G_{x0}}}
\def\Nx0{{N_{x0}}}
\def\Nxyz{{N_{xyz}}}
\def\N0{{N_{0}}}
\def\Qx0{{Q_{x0}}}
\def\Ciso{{2^{1+24}_+.\Co_1}}
\def\Noiso{{2^{2+11+22}.(M_{24} \times S_3)}}

\begin{document} {\large}

\maketitle

%\centerline{\large Preliminary version}

%{\noindent \bf Key Words:}
%Monster group, finite simple groups, group representation

% \vspace{1ex}

%\noindent MSC2020-Mathematics Subject Classification (2020):
%20C34, 20D08, 20C11

\begin{center}
\textit{In memory of John Conway and Simon Norton, whose work continues to inspire.}
\end{center}
\bigskip % Adds some vertical space before the abstract

\begin{abstract}
\small 

We determine the order of the largest of the twenty-six sporadic simple groups known as the Monster, using a straightforward computational approach.
The Monster is here defined as a subgroup of the symmetry group of the $196884$-dimensional Griess algebra generated by
a group of type $2^{1+24}_+.\Co_1$ and an additional triality automorphism.
Our approach is based on counting arguments for certain idempotents of the Griess algebra called {\em axes}.
Our proof is self-contained, requiring only established properties of the Conway group as the automorphism group of the Leech lattice,
and some of its subgroups. 

Although our approach is conceptually simple, it requires extensive calculation inside a $196884$-dimensional matrix group 
that current computer algebra systems cannot easily handle directly. 
Instead, we use the software package {\em mmgroup}, developed by the second author, which supports fast calculations inside the Monster.
To our knowledge, this paper contains the first self-contained computation of the order of the Monster.

The Monster also acts on the Moonshine module $\V$, which is a vertex operator algebra of central charge $c=24$.
We provide a new proof that the Monster is the {\em full automorphism group} of the Griess algebra and of the Moonshine module using Borcherds' proof of the Monstrous Moonshine conjectures.
In addition, we show that the Monster has exactly two conjugacy classes of involutions.
The order of the Baby Monster, the second largest of the sporadic simple groups, is also determined.

\end{abstract}

%%%%%%%%%%%%%%%%%%%%%%%%%%%%%%%%%%%%%%%%%%%%%%%%%%%%%%%%%%%%%%%%%%%%%%%%%%%%%%

	\section{Introduction}
\label{section:Introduction}

In this paper, we will determine the order of the Monster finite group $\M$
assuming only established properties of the Leech lattice and its automorphism group, and using extensive computations 
performed with the software package {\em mmgroup}~\cite{mmgroup2020} developed by the second author.

\smallskip
The Monster is the largest of the twenty-six sporadic finite simple groups. 
Its existence was predicted independently at the end of 1973 by Bernd Fischer 
and Robert Griess. Shortly afterwards, its presumptive order was computed to be
$$  2^{46}\cdot 3^{20}\cdot 5^9\cdot 7^6 \cdot  11^2 \cdot  13^3 \cdot 17 
\cdot 19 \cdot 23\cdot  29 \cdot  31\cdot   41 \cdot  47 \cdot  59 \cdot  71 \qquad\qquad \qquad\qquad\qquad $$
$$\qquad  =\  808,017,424,794,512,875,886,459,904,961,710,757,005,754,368,000,000,000. $$

\medskip

The Monster was first constructed by Griess in 1982~\cite{Griess:Friendly:Giant} as a group of automorphisms
of a $196884$-dimensional non-associative algebra, called the {\em Griess algebra} and denoted by~$\Griess$.
It was later discovered that the Griess algebra is part of a larger infinite-dimensional structure
called the {\em Moonshine module} and denoted by $\V$~\cite{Bo-ur,FLM}.
The Moonshine module is an extremal vertex operator algebra of central charge~$24$ with the Monster 
as its symmetry group.
Vertex operator algebras are a mathematical formalization of two-dimensional conformal
field theories first studied in mathematical physics.

The automorphism group of the Griess algebra --- 
or, alternatively, the Moonshine module --- 
is known to contain a visible subgroup $\C$ of type $\Ciso$ (the centralizer of an involution $x_{-1}$) and an additional {\em triality automorphism} $\tau$ of order~$3$.
The Monster is defined as the subgroup generated by these:
$$ \M=\langle \C, \tau \rangle \subseteq  \Aut(\Griess).$$

Several group-theoretical properties of the Monster, like its finiteness, can be understood in
terms of vertex operator algebras. 
Carnahan has recently provided upper and lower bounds for the order using
basic group theory, the theory of orbifolds of vertex operator algebras, and Borcherds' 
proof of the Monstrous Moonshine conjectures~\cite{Carnahan:2023}.
However, the order of the Monster has never been determined directly by using only properties of 
the Griess algebra or the Moonshine module alone.

\medskip

The only published proof for the order given above of the Monster can be found in 
\cite{griess1989uniqueness}, Corollary~3.7.3. 
One starts by identifying the double cover $2.B$ of the Baby Monster as the centralizer of an involution in $\M$.
Then the order of the Monster is obtained by enumerating the nine double cosets in
$2.B \backslash \M / 2.B$, and by computing the
number of cosets in $2.B \backslash \M$ contained 
in each of these double cosets. 
In \cite{griess1989uniqueness}, properties of the Baby Monster $B$ as well 
as of two Conway groups $\Co_1$ and $\Co_2$, the Thompson group ${\rm Th}$, the Harada-Norton group ${\rm HN}$,
and the two Fischer groups ${\rm Fi}_{22}$ and ${\rm Fi}_{23}$ were used. In particular, their orders are assumed to be known in advance.
The obvious disadvantage is that all of these groups have to be studied beforehand, either theoretically or computationally.

\bigskip

%%%%%%%%%%%%%%%%%%%%%%%%%%%%%
\paragraph{A strategy for computing the order of the Monster.}

Our approach in this paper is conceptually much simpler. We consider the permutation action of 
$\M$ on certain idempotents in $\Griess$ called {\em axes} --- or {\em Ising vectors} if considered in $\V$. 
%These are the elements in $\Griess$ that generate a Virasoro vertex operator subalgebra of central charge $c=\frac{1}{2}$ in $\V$.
Associated with each axis is a 2A involution in $\M$.

Denote by $\M_{u}$ the stabilizer of an axis $u$ in $\M$.
Let $X^+$ be the $\M$-orbit of an explicitly given 
axis $v^+$ and
let $X^-$ be the $\M_{v^+}$-orbit of an explicitly given axis $v^-$ in $X^+$ orthogonal to $v^+$. 
We will determine the sizes of $X^+$ and $X^-$. 
Using the fact that the common stabilizer $\M_{v^+} \cap \M_{v^-}$ 
of $v^+$ and $v^-$ is contained in $\C$, one shows $\M_{v^+} \cap \M_{v^-}\cong 2^{1+23}.\Co_2$.
It follows that 
$$|\M| = |X^+|\cdot |\M_{v^+}| = |X^+| \cdot |X^-| \cdot |2^{1+23}.\Co_2|.$$

\smallskip

The basic strategy to determine the size of $X^+$ and similarly $X^-$ is as follows.
The decomposition of $X^+$ into twelve $\C$-orbits 
and the structure of their stabilizers has been given in~\cite{Norton98} without proof. 
One can decompose these orbits further into orbits under the group 
$\Nx0\cong 2^{1+24}_+.(2^{11}{:}M_{24})\cong 2^{2+11+22}.(M_{24} \times 2)< \C$ by considering the action of stabilizers of the twelve 
$\C$-orbits on  $\Lambda/2\Lambda \cong 2^{24}$, where $\Lambda$ denotes the Leech lattice.
Using the software package mmgroup, we provide 
explicit representatives of axes in the Griess algebra for each of these $\Nx0$-orbits and 
explicit generators for the stabilizers of the twelve $\C$-orbits.

The group $\Nx0$ is a subgroup of index~$3$ in the group $\N0 =\langle \Nx0, \tau \rangle$ of structure $2^{2+11+22}.(M_{24} \times S_3)$.
Let  $\Nxyz$ be the subgroup of $\Nx0$ of index 2 and structure $2^{2+11+22}.M_{24}$.
Then $\Nxyz$ is normal in both $\Nx0$ and $\N0$; but $\Nx0$ is not normal
in $\N0$, see e.g.~\cite{Conway:Construct:Monster}.

Since $\Nxyz \vartriangleleft \N0$, the factor group $N_0/\Nxyz\cong S_3$
permutes the $\Nxyz$-orbits contained in an $\N0$-orbit. This fact
will help us to show that there are $251$ $\Nx0$-orbits on the axes fusing
into $123$ $\N0$-orbits and twelve $\C$-orbits.
Those orbits fuse into a single orbit $X^+$ of axes closed under the action
of both $\C$ and the triality automorphism~$\tau$.

A counting argument allows us to determine the sizes of the twelve $\C$-orbits and of~$X^+$.
We also show that our explicit elements inside the stabilizers of each of the twelve $\C$-orbits 
actually generate the whole stabilizer. 
\medskip

An analogous approach is then taken for the action of $\M_{v^+}$ on $X^-$. 
There are ten $(\M_{v^+}\cap \C)$-orbits in $X^-$ and 
the triality automorphism stabilizes $v^+$ and fuses the ten orbits.
\smallskip

These calculations depend extensively on calculations in mmgroup describing the action of $\C$
and the triality automorphism on the Griess algebra $\Griess$. We also use
calculations in $\C$ done with the software package Magma~\cite{bosma1997magma}. The only information from
other sporadic groups used is the knowledge of the Conway groups $\Co_1$ and $\Co_2$ obtained from
the automorphism group of the Leech lattice, and of the Mathieu group $M_{24}$ as the automorphism group
of the binary Golay code of length~$24$.

\medskip

Our approach is analogous to Conway's original approach to the construction of $\Co_0$ and the computation of 
its order~\cite{Conway-conway}, and it generalizes this method to the third generation of the happy family~\cite{Griess:Friendly:Giant}.

It may be possible to simplify our approach further by using additional 
structure-theoretical results about vertex operator algebras. For example, if it can be shown directly 
that the twelve $\C$-orbits contain all axes of $V_2^\natural$, one can use arguments similar 
to those in~\cite{Conway-conway} that these orbits must fuse under the additional triality automorphism.

\smallskip

Combining our result with the result of Carnahan~\cite{Bo-ur,Carnahan:2023}, we also show that the Monster
is the full automorphism group of $\Griess$ and $\V$. The previous proofs of this result given 
by Tits~\cite{ti-monster,ti-monster2} all depend on certain group-theoretic characterizations of the Monster.

\bigskip

%%%%%%%%%%%%%%%%%%%%%%%%%%%%%%%%%%%%%

\paragraph{Outline.} 
The paper is organized as follows. In the next section, we review basic 
properties of the Griess algebra and of the Moonshine module and 
discuss the visible parts of their automorphism groups.
In Section 3, we describe the $251$ different $\Nx0$-orbits on the axes
and how these orbits fuse under $\C$ and $\N0$. 
In particular, we obtain the number $|X^+|$ of axes in the Griess algebra. 
In Section 4, we repeat a similar analysis for the stabilizers $\M_{v^+}$
of $v^+$ acting on $X^-$. 
In particular, we obtain the number $|X^-|$ of axes in the
$\M_{v^+}$-orbit of an axis orthogonal to a given axis $v^+$.
In the final Section~\ref{section:order:monster}, we provide several applications, including the 
computation of the order of the Monster and the Baby Monster.
The appendices provide details of the mmgroup calculations. In Appendix~C, the orbits under
the action of the maximal subgroup $\N0$ on the axes are listed.

\smallskip

As a supplement to this paper, we provide the programs~\cite{accompanying} % This may be ok, but maybe only when the  code is in a different location than the paper.
which allow to verify all computational results obtained with mmgroup.

%%%%%%%%%%%%%%%%%%%%%%%%%%%%%%%%%%%%%%%%%%%%%%%%%%%%%%%%%%%%%%%%%%%%%%%%%%%%%%
\section{Construction of the Monster}

The Monster can be defined either as a group of automorphisms of the $196884$-dimensional 
Griess algebra $\Griess$ as in~\cite{Griess:Friendly:Giant, Conway:Construct:Monster} or as a group of automorphisms of
the Moonshine module vertex operator algebra $\V$ of central charge $24$ as in~\cite{Bo-ur,FLM}.
In both approaches, the Monster is defined as the subgroup generated by two ``visible'' subgroups $\C$ and $\N0$.

\medskip

In the first subsection, we will review the construction of the Monster and the Griess 
algebra as in Conway's paper~\cite{Conway:Construct:Monster}.
%Our calculation of the order of the Monster uses the software package {\em mmgroup} of
%the second author which is based on~\cite{Conway:Construct:Monster}.
We also provide proofs for certain properties of involutions and their centralizers in the Monster.
Although these properties can also be deduced more abstractly from the theory of vertex operator algebras,
we prefer to provide full proofs to make our computations of the order of the Monster more self-contained.

\medskip

In the second subsection, we will review the construction of the Monster in terms of the
Moonshine module. We need this description to be able to provide a new proof of the result of Tits
that the Monster is the full automorphism group of the Griess algebra and of the Moonshine
module. We will also explain how the above-mentioned properties of centralizers of involutions
can be shown using vertex operator algebras.

\medskip

In the final subsection, we identify the automorphism group of the Griess algebra with the automorphism group of the Moonshine module.

%%%%%%%%%%%%%%%%%%%%%%%%%%%%%%

\subsection{The Monster and the Griess algebra}
\label{Monster:Griess}

We first define a $196884$-dimensional common rational representation
$\Griess$ of two groups $\C$
of structure $\Ciso$ and $\N0$ of structure $\Noiso$ with matrix entries 
in~$\setZ[\frac{1}{2}]$, as in \cite{Conway:Construct:Monster}.
The Monster $\M$ is defined to be the group generated by $\C$ and $\N0$. 
The group $\Nx0 = \N0 \cap  \C$ has structure $2^{1+24+11}.M_{24}$
and index~$3$ in $\N0$, see~\cite{Conway:Construct:Monster}.
We write $\Qx0$ for the normal subgroup  of $\C$ of structure $2^{1+24}_+$,
which is extraspecial of plus type, and $x_{-1}$ for the unique central
involution in $\Qx0$. We write $\Lambda$ for the Leech lattice and
$\Lambda / 2 \Lambda$ for the Leech lattice mod 2. Then
$\Aut(\Lambda / 2 \Lambda)$ is the simple group $\Co_1$; and
$\Aut(\Lambda)$ is the group $\Co_0$ of structure $2.\Co_1$.  
Since $\Co_1$ is simple, $x_{-1}$ is also the
unique central involution in $\C$.

\smallskip

There is a natural homomorphism from $\Qx0 \cong 2_+^{1+24}$ to
$\Lambda / 2 \Lambda$ with kernel $\{1,\, x_{-1}\}$.
The {\em type} of a vector  in $\Lambda / 2 \Lambda$ is the halved
squared norm of its shortest preimage in $\Lambda$; possible types
are $0$, $2$, $3$, and $4$. The type of an element
of $\Qx0$  is the type of its image in $\Lambda / 2 \Lambda$.
The  group  $\Co_1$ is transitive on the vectors of any given type
in $\Lambda / 2 \Lambda$.
For $r \in \Qx0$ we write $x_r$  when we consider $r$
as an element of $\C$ or $\M$ via the injective
mappings $\Qx0 \hookrightarrow \C \hookrightarrow \M$,
as in~\cite{Conway:Construct:Monster,Seysen20}; and we write
$\lambda_r$ for the image of $r$ in  $\Lambda / 2 \Lambda$.

\smallskip

Elements of $\Lambda / 2 \Lambda$ or  $\Qx0$ of  type~$2$  are
called {\em short}. For a short $r \in \Qx0$ we also write $\lambda_r$
for a shortest vector in the Leech lattice that maps to the same
element of $\Lambda / 2\Lambda$ as $r$; so $\lambda_r$
has squared norm~$4$ in $\Lambda$.
By standard properties of the Leech lattice, $\lambda_r$ as an element of
$\Lambda$ is defined up to sign only; but $\lambda_r \otimes \lambda_r$
is uniquely defined as an element of the symmetric tensor square
$S^2(\Lambda \otimes_\Z \Q)$.

\smallskip

Let the element $x_\Omega$ of  $\Qx0$ be as in
\cite{Conway:Construct:Monster} and \cite{Seysen20}.
Then $x_\Omega$ is of type~$4$; and the image
$\lambda_\Omega$ of  $x_\Omega$ in $\Lambda / 2 \Lambda$
is equal to $v + 2 \Lambda$ for any primitive $v \in \Lambda$
proportional to a standard basis vector of $\Lambda$.
Let $x_{\pm \Omega}$ be the two preimages of $\lambda_\Omega$
in $\Qx0$.
As in \cite{Conway:Construct:Monster},
we embed $\N0$ into $\M$ so that the characteristic subgroup of
structure $2^2$ of $\N0$  is equal to $\{x_{\pm 1},\, x_{\pm \Omega}\}$.

\medskip

We take for granted all theorems proved in   
\cite{Conway:Construct:Monster}, §§ 1--12, and Appendices \hbox{1--7}.
Apart from this, we use well-known facts about the
Mathieu group $M_{24}$ and the Leech lattice and
its automorphism group as stated e.g.\ in
\cite{Conway-SPLG}, Ch.~4.11, and 10--12, or
\cite{citeulike:Monster:Majorana}, Ch.~1.
We also assume that character information for the group
$\Co_1$ and for its sections $\Co_2$ and 
$M_{24}$ from the ATLAS \cite{Atlas} is available. But we do
not use the character information for the Monster or the Baby
Monster in \cite{Atlas}, since it is difficult to understand 
how this information is obtained without assuming the existence
of the Monster.

\smallskip

Our proof relies on explicit computations in the representation
$\Griess$ of $\M$, which is called  $196884_x$ in 
\cite{Conway:Construct:Monster}. 
For computations in $\Griess$ we use the mmgroup package~\cite{mmgroup2020, mmgroup_doc},
which supports computations in~$\Griess$. We also adopt the notation
from~\cite{Seysen20}, where $\Griess$ is called $\rho$, and
explicit generating systems for $\C$, $\N0$, and $\M$ are given. 
Computations in $\Griess$ are done modulo~$15$.
This is possible because the denominators of the entries of the matrices
representing elements of $\M$ are powers of~$2$, see e.g.~\cite{Conway:Construct:Monster, Seysen20}. 
For a detailed discussion of why computation modulo~$15$ is preferred, we refer to~\cite{Seysen22}.

\medskip

%%%%%%%%%%%%%%%%%%%%%%%%%%%%%%%%%%%%%

\paragraph{The Monster representation.}
To construct the representation $\Griess$ of the Monster 
we start with a representation of the group $\C$
on the rational vector space
\[
\Griess = 300_x \oplus 98280_x \oplus 98304_x \,, \quad \mbox{where}
\]

\vspace{-6.5mm}
\begin{align*}
300_x & \mbox{ is the symmetric tensor square 
	$S^2(24_x)$ of } \, 24_x, \\
24_x & \mbox{ is the representation } \Lambda \otimes_\Z \Q \mbox{ of }  
  \Co_0  \;
  \mbox{as the automorphism group of } \Lambda  ,\\
98280_x & \mbox{ is the monomial representation of } \C\; 
\mbox{ acting on short elements of } \Qx0 , \\
98304_x & \mbox{ is the tensor product } \; 4096_x \otimes 24_x \,
\mbox{ with } 4096_x 
\mbox{ as in \cite{Conway:Construct:Monster,Seysen20}.}
\end{align*}

A basis vector of $98280_x$ is denoted by $X_r$, where
$r$ is a short element of $\Qx0$.
%We write $x_r$ for $r$ when we consider $r$
%as an element of $\C$ or $\M$ via the injective
%mappings $\Qx0 \hookrightarrow \C \hookrightarrow \M$,
%as in   \cite{Conway:Construct:Monster,Seysen20}.
We identify $X_{r \cdot x_{-1}}$ with
$-X_{r}$.   As in
\cite{Conway:Construct:Monster}, we write $\widetilde{\Omega}$
for the set of size~$24$ on which the Mathieu group $M_{24}$ acts,
and $ \mathcal{P}$ for the {\em Parker loop}.
There the natural basis vectors of $24_x$ are denoted by
$i_1$, $j_1$, $\ldots$, for $i$, $j$, $\ldots \in \widetilde{\Omega}$.
So we have the following basis vectors of $\Griess$:
\begin{align*}
\mbox{for } \hphantom{00}300_x: & \quad
(ii)_1 := i_1 \otimes i_1  \, , \qquad  
( i \in \tilde{\Omega}),  \quad & \mbox{of norm $1$}, \\    
& \quad
(ij)_1 := i_1 \otimes j_1 + j_1 \otimes i_1 \, , 
 \quad ( i,\, j \in \widetilde{\Omega},\   i \neq j)
\, ,   &   \mbox{of squared norm $2$}, \\    
\mbox{for } 98280_x: & \quad X_r \, , \quad
 (r \in \Qx0 \; \;  \mbox{short})
\, ,   &   \mbox{of norm 1}, \\ 
\mbox{for } 98304_x: & \quad
d_1^\pm \otimes i_1, \quad
(d \in \mathcal{P},\  i \in \widetilde{\Omega}  ) 
\, ,   &   \mbox{of norm 1}. 
\end{align*}
For the definitions of the basis vectors of $98304_x$
we refer to \cite{Seysen20}. These basis vectors are
orthogonal, except when equal or opposite. 
The basis vectors of $ 300_x \oplus 98280_x$ in
\cite{Conway:Construct:Monster} agree with those in
 \cite{Seysen20}; the other basis vectors of $\Griess$ agree
up to sign.

The representation $300_x$ of $\C$ has a natural 
interpretation as the space of symmetric matrices over $24_x$.
It has a decomposition $300_x = 1_x \oplus 299_x$, where $1_x$
is the space spanned by the unit matrix, and $299_x$ is the space
of the symmetric matrices of trace zero.  
We write $1_\Griess$ for the element 
of $300_x$ corresponding to the unit matrix divided by four.

To obtain a representation $\Griess$ of the Monster $\M$, we
also have to describe the action of a {\em triality} element $\tau$
on $\Griess$, such that $\N0$ is generated by $\Nx0$ and $\tau$. 
For details we refer to \cite{Conway:Construct:Monster,Seysen20}.
Then $\M$ fixes $1_\Griess$ in $\Griess$; so we have a
decomposition
\[
   \Griess =  1_x \oplus 196883_x, \quad \mbox{with} \quad
      196883_x =  299_x \oplus 98280_x \oplus 98304_x \; ,
\] 
where $1_x = {\rm Span}(1_\Griess)$ and $196883_x$ are representations of $\M$.

\smallskip
We write $(a,b)$ for the Monster-invariant scalar product of 
two elements $a$, $b$ of $\Griess$ as defined in~\cite{Conway:Construct:Monster}.
This is to be distinguished from the scalar product 
$\langle ., . \rangle$ in $\Lambda/2 \Lambda$. 
We also need an algebra structure on $\Griess$ called the
{\em Griess algebra}, which is invariant under the action of $\M$. 
We adopt the definition of the Griess algebra from \cite{Seysen20},
which is compatible with the corresponding definition in   
\cite{Conway:Construct:Monster}, §~12. We denote the algebra
product of $a$ and $b$ by $a * b$. That
product is commutative, but not associative. We write
$(a,b,c)$ for $(a * b, c)$. Then $(a,b,c)$ is a symmetric
trilinear form on $\Griess$.  The element
$1_\Griess$ of squared norm $\frac{3}{2}$ acts as the identity of the
Griess algebra. We write $\ad a$ for the
endomorphism $v \mapsto a * v$ of the vector space $\Griess$.

\medskip

%%%%%%%%%%%%%%%%%%%%%%%%%%%%%%%%%%%%%

\paragraph{Involutions and axes.}

By standard properties of $\Qx0$ and the Leech lattice, all short 
elements in $\Qx0$ are in the same conjugacy class in $\C$.

\begin{Lemma}
The centralizer of $x_r$ (with $r \in \Qx0$ short) in $\C$
fixes precisely the subspace of $\Griess$ spanned by
$1_\Griess$, $\lambda_r \otimes \lambda_r$, and~$X_r$.  
\end{Lemma}

\fatline{Proof:}
Write $C_r$ for the centralizer of $x_r$ in $\C$ and $\Griess_r$
for the subspace of $\Griess$ fixed by $C_r$. Obviously,
$1_\Griess$, $\lambda_r \otimes \lambda_r$, and $X_r$ are in $\Griess_r$.
The element $x_{-1}$ of $C_r$  negates $98304_x$; hence
$\Griess_r \subset  300_x \oplus 98280_x$. Assume that $v \in  \Griess_r$
has a nonzero co-ordinate at the basis vector $X_s$, $s \neq r$,
$r \in \Qx0$ short. Then we can find a $t \in \Qx0$ 
(with  image $\lambda_t$ in $\Lambda/2 \Lambda$) such that
$\langle \lambda_t , \lambda_r \rangle = 0$,
$\langle \lambda_t, \lambda_s \rangle= 1 \pmod{2}$.
Since $\Lambda$ is transitive on pairs of short vectors with a given
scalar product, it suffices to find suitable elements $t$ of $\Qx0$  for
one pair $\lambda_r, \lambda_s$ with each possible scalar product
$\langle \lambda_r, \lambda_s \rangle$. Possible scalar products are
$0$, $\pm 1$, $\pm 2$, and $\pm 4$; so we can find suitable
$t \in \Qx0$ computationally by checking random elements of $\Qx0$.
Any such $t$ fixes $X_r$ (so it also centralizes $x_r$)
and negates $X_s$. So the co-ordinate
of $v$ at the basis vector $X_s$ must be zero, a contradiction. Thus
$C_r$ fixes at most $300_x \oplus  \langle X_r \rangle$.

The centralizer of $\lambda_r$ in $\Co_1 = \C/ \Qx0$
is the simple group $\Co_2$. The preimage of 
$\Co_2$ in $\C$ either fixes or negates $\lambda_r$.
Thus $300_x$ is also a representation of   
$\Co_2$. Character calculations in $\mbox{Co}_1$ and
$\Co_2$ show that the character of $300_x$ as a character
of $\Co_2$ is a sum 
$\chi_1 \oplus \chi_1 \oplus \chi_{23} \oplus \chi_{275}$
of the irreducible characters $\chi_1$, $\chi_{23}$, and $\chi_{275}$
of $\Co_2$ of dimensions $1$, $23$, and $275$. Thus the subspace
of $300_x$ fixed by $C_r$ 
%fixing $\Griess_r$ 
has dimension~$2$. \hfill \proofend

\smallskip
 
\begin{Corollary}
\label{cor:cent_x_r}	
The centralizer of $x_r$ (with $r \in \Qx0$ short) in $\M$
fixes a one- or two-dimensional subspace of $\Griess$ spanned by
$1_\Griess$ and, possibly, another vector $\lambda_r \otimes \lambda_r - 4 X_r$.
\end{Corollary}
 
\fatline{Proof:} The Monster $\M$ fixes $1_\Griess$. 
Since all involutions $x_r$ in $\Qx0$ with $r$ short are conjugate in
$\C$, we may assume $x_r = x_{ij}$ in the notation 
of~\cite{Conway:Construct:Monster, Seysen22}.
The triality element $\tau$ centralizes $ x_{ij}$.
In that notation we have
$\lambda_r \otimes \lambda_r = 2(ii)_1 + 2(jj)_1 - 2 (ij)_1$.
Let
$V = {\rm Span}( 1_\Griess,\, \lambda_r \otimes \lambda_r,\, X_r )
\subset \Griess$.
From the action of $\tau$ described in \cite{Seysen22},
Section~8.1, we see that the subspace of $V$ fixed by
$\tau$ is just the space ${\rm Span}(1_\Griess,\, \lambda_r \otimes \lambda_r - 4 X_r)$.
\hfill \quad \proofend

\smallskip

We want to show that the subspace fixed by the centralizer of $x_r$ 
in the Corollary is in fact $2$-dimensional. 

\begin{Definition}
   The {\em axis} of $x_r$ (with $r \in \Qx0$ short)
	is the vector $ \ax(x_r)= \frac{1}{2}\lambda_r \otimes \lambda_r - 2 X_r$
	in $\Griess$.
\end{Definition}

\begin{Lemma}
\label{Lemma:pre:central:2A:axis}
Let $v = \ax(x_r)$ (with $r \in \Qx0$ short).
Then $v * v = 16\,v$ and $(v,v) = 8$.
The endomorphism $\ad v$ of $\Griess$
has eigenvalues $16$, $0$, $4$, and~$\frac{1}{2}$ with
eigenspaces of dimensions $1$, $96256$, $4371$, and $96256$, respectively.
The involution $x_r$ negates the eigenspace of the
eigenvalue~$\frac{1}{2}$ and fixes the other eigenspaces
of  $\ad v$. 
\end{Lemma}

\fatline{Proof:} The proof is similar to the proof in~\cite{Conway:Construct:Monster}, Appendix~6. 
Since all involutions $x_r$ in $\Qx0$ with $r$ short are conjugate in
$\C$, we may assume $x_r = x_{ij}$,  $i$, $j \in \widetilde{\Omega}$ in the
notation  of~\cite{Conway:Construct:Monster}. Then $x_{ij}$ either fixes
or negates a basis vector of $\Griess$. A simple calculation shows $(v,v) = 8$
for $v = \ax(x_{ij})$.

Let $V_X$, $V_Y$, and $V_Z$ be the subspaces of $\Griess$ defined in
\cite{Conway:Construct:Monster},~§11. Then $V_X$ is the space
spanned by the basis vectors $X_s$ of $98280_x$ with
$  \langle\lambda_\Omega, \lambda_s \rangle = 1  \pmod 2$. We also have
$V_Y =  (V_X)^\tau$, $V_Z = (V_X)^{\tau^{-1}}$, and
$98304_x =  V_Y \oplus V_Z$, see \cite{Conway:Construct:Monster},~§11.
Here $\tau$ is the triality element in $\M$ as above.
We define $V_0$ to be the subspace of $300_x \oplus 98280_x$ spanned
by $300_x$ and the basis vectors of $98280_x$ with
$  \langle\lambda_\Omega, \lambda_s \rangle = 0  \pmod 2$. So we obtain
the following decomposition of $\Griess$:

$$\begin{array}{|c|c|c|c|c|}
300_x &  \multicolumn{2}{c|}{ 98280_x}  &
     \multicolumn{2}{c|}{ 98304_x } \\
\hline
   \multicolumn{2}{|c|}{ \phantom{mm} V_0 \phantom{mm}}  & \phantom{mm} V_X \phantom{mm} &
    \phantom{mm} V_Y \phantom{mm} & \phantom{mm} V_Z \phantom{mm}
\end{array} \quad .$$

Direct calculation shows that $\ad v$ has the following eigenvalues and eigenvectors
where  $r$, $s$, and $t \in \Qx0$  are short and disjoint:
\newcommand{\nEvenOdd}[2]{{#1 \small{\mbox{ even }} + #2 \small{\mbox{ odd}}}}
\renewcommand{\arraystretch}{1.2}
$$\begin{array}{r|ll|rr}
\multicolumn{1}{l|}{\mbox{eigen-}} & \mbox{eigenvectors}  & & 
  \multicolumn{2}{l}{\mbox{dimension of subspace}} \\ 
\multicolumn{1}{l|}{\mbox{value}} &         & 
     & \mbox{in } V_0 & \mbox{in } V_X \\ \hline
16 & \textstyle  \frac{1}{2} \lambda_r \otimes \lambda_r - 2 X_r  & & 1  \\
0  & {\textstyle \frac{1}{2}} \lambda_r \otimes \lambda_r + 2 X_r  & & 1 \\
4  &  \lambda_r \otimes \lambda_s + \lambda_s \otimes \lambda_r &\mbox{if } 
 \langle \lambda_r, \lambda_s \rangle = 0  &  23\\
0  & \lambda_s \otimes \lambda_t + \lambda_t \otimes \lambda_s 
 &  \mbox{if }   \langle\lambda_r, \lambda_s \rangle 
   =  \langle\lambda_r, \lambda_t \rangle = 0&  \tiny{{24}\choose{2}} = 276  \\
0 & X_s + X_t & \mbox{if } rst = 1 & 1276 & 1024 \\
4 & X_s - X_t & \mbox{if } rst = 1 &  1276 & 1024 \\
{\textstyle \frac{1}{2}} &   X_s  & \mbox{if } \langle\lambda_r, \lambda_s \rangle = \pm 1
 & 22528 & 24576 \\
0 & X_s & \mbox{if }  \langle\lambda_r, \lambda_s \rangle = 0 & 24047 & 22528
\end{array}$$
The scalar product $\langle\lambda_r, \lambda_s \rangle$ is
$0$, $\pm 1$, $\pm 2$, or $\pm 4$. It is $\pm 4$ in case $r = \pm s$.
If it is $\pm 2$ then there is a short $t$ in $\Qx0$ with
$rst = 1$. Thus any vector $w$ in  $300_x \oplus \break {\rm span}( \{X_s \mid s \in \Qx0 \mbox{ short}\})$
is a linear combination of the eigenspaces of $\ad v$ stated
in the table above. So these vectors $w$  span $300_x \oplus 98280_x$.

In the subspace $300_x \oplus 98280_x$ of $\Griess$ the involution
$x_r$ negates precisely the subspace spanned by basis vectors $X_s$ of
$98280_x$ with $\langle \lambda_r, \lambda_s \rangle = 1 \pmod{2}$.
From the table above we see that the negated basis vectors
are in the eigenspace of the eigenvalue $\frac{1}{2}$ of
$\ad v$, and that the fixed basis vectors are linear
combinations of the other eigenspaces.
Since $\tau$ centralizes
$x_{ij}$, and the basis vectors of $98304_x$ are images of the basis vectors in  $98280_x$
under $\tau^{\pm1}$, this is true for all basis vectors of $\Griess$.  

So the dimensions of the intersections of the eigenspaces of $\ad v$
with $V_0 \oplus V_X$ can be obtained from the table above by summing up
all entries of the last column for each eigenvalue. 
Since $x_r$  commutes with $\tau$, the intersections of an eigenspace of
$\ad v$ with $V_X$, $V_Y$, and $V_Z$ all have the same dimension. So we have
to count the entries marked as `in $V_X$' in the last column three times
in order to obtain the dimension of an eigenspace of $\ad v$. The entries for 
the basis vectors $X_s$ have been computed with the Python script
{\tt eigenspaces\_axis.py} in the accompanying program code; but they
can also be obtained from standard properties of the Leech lattice. 
\hfill \proofend

\smallskip

\fatline{Remark:}
As a byproduct of the proof of Lemma~\ref{Lemma:pre:central:2A:axis}
we have obtained the dimensions of the eigenspaces of $\ad v$. These
dimensions have also been computed in
\cite{citeulike:Monster:Majorana}, Ch.~8.4. We have scaled the largest
eigenvalue of  $\ad v$ to $16$ as in \cite{Seysen22}. That eigenvalue
is scaled to $64$ in \cite{Conway:Construct:Monster}, and to $1$ in
\cite{citeulike:Monster:Majorana}.

\smallskip

\begin{Lemma}
\label{lemma:def:axis}	
The centralizer $C_{\M}(x_r)$ of the involution $x_r$ where $r \in \Qx0$ is short 
fixes a two-dimensional subspace of $\Griess$ spanned by
$1_\Griess$ and the vector $\ax(x_r)$. The vector $\ax(x_r)$ is the unique vector $v$ 
in that subspace satisfying $v * v = 16 \cdot v$  and $(v,v) = 8$.	
\end{Lemma}

\fatline{Proof:} Let $V^-$ be the $-1$ eigenspace of $x_r$. 
For any  $w \in V^-$ of norm~$1$ one has 
$(w*w,1_\Griess)=(w,w*1_\Griess)=(w,w)=1$. By
Lemma~\ref{Lemma:pre:central:2A:axis} we have
$(w*w,\ax(x_r)) = (w, w * \ax(x_r)) = \frac{1}{2} (w,w) = \frac{1}{2}$.  
Averaging $w*w$ over all vectors $w$ in $V^-$ of norm~$1$  yields a
linear form on $\Griess$ corresponding to a 
vector $u \in \Griess$ invariant under $C_{\M}(x_r)$  with
$(u,1_\Griess) = 1$ and $(u, \ax(x_r)) = \frac{1}{2}$. 
On the other hand we have $(1_\Griess, 1_\Griess) = \frac{3}{2}$ and 
$(1_\Griess, \ax(x_r)) = \frac{1}{2}$. So $u$ and $1_\Griess$ are
not proportional; hence the subspace $V_2$ of $\Griess$ fixed by
$C_{\M}(x_r)$ has dimension at least~$2$. 
By Corollary~\ref{cor:cent_x_r} that dimension must be equal to~$2$.

Put $v' = \ax(x_r) / 16$. Then  $1_\Griess$, $v'$, and~$1_\Griess - v'$
are idempotents in $V_2$ with norms $\frac{3}{2}$, $\frac{1}{32}$,
and~$\frac{3}{2}$, respectively.
Since $v' * (1_\Griess - v') = 0$, there are no further idempotents
in $V_2$. Thus  $v = \ax(x_r)$ is the unique vector in $V_2$ with
 $v * v = 16 \cdot v$  and $(v,v) = 8$. 
\hfill  \proofend

\smallskip

We remark that the proof of the corresponding statement in
\cite{Conway:Construct:Monster}, §~14 uses the character
table of the Monster.
We may now define:
\begin{Definition}\rm
\begin{itemize}
\item[(a)]
	A 2A {\em involution} of the Monster is an element conjugate in $\M$ to~$x_r$ for a short element $r\in \Qx0$.
\item[(b)]
   The {\em axis} $\ax(t)$ of a 2A involution $t$ is the unique vector
    $v \in \Griess$ that is fixed by $C_{\M}(t)$ and
    satisfies $v * v = 16 \cdot v$  and $(v,v) = 8$.
\item[(c)]
    The set of all such vectors $\ax(t)$ forms the set of {\em axes} of $\Griess$.
\end{itemize}
\end{Definition}

The axis of a 2A involution is well-defined by Lemma~\ref{lemma:def:axis}.
By Lemma~\ref{Lemma:pre:central:2A:axis}
a 2A involution $t$ can be recovered from its axis $\ax(t)$
as the map that negates the $\frac{1}{2}$ eigenspace of
$\ad(\!\ax(t))$ and fixes the remaining eigenspaces.
This establishes a one-to-one correspondence between the 2A
involutions and their axes.
\smallskip

Recall that the involutions $x_r$ for short 
elements $r\in \Qx0$  are conjugate in $\C$ and hence in $\M$.
Thus all 2A involutions are conjugate in $\M$.
However,  we cannot exclude yet that the conjugacy
class of a 2A involution in $\Aut(\Griess)$ is greater than
the class of 2A involutions in $\M$.

It is well known that the quotient of the centralizer of a 2A involution
$t$ in $\M$ by the subgroup $\langle t \rangle$ is a simple group $B$,
known as the {\em Baby Monster}, the second largest of the sporadic simple groups; cf.~\cite{Griess:Friendly:Giant}, Lemmas~13.2 and~13.3.
In this paper, we define $B$ as this quotient, without assuming \emph{a priori} that $B$ is simple.

\medskip

As an application, we show:
\begin{Theorem}
The automorphism group of the Griess algebra and hence the Monster is finite.
\end{Theorem}

\fatline{Proof:} The same proof as in  \cite{Conway:Construct:Monster}, §~13, and
Appendix~6 works. The images
of an axis $\ax(x_r)$ in $\M $ span $300_x \oplus 98280_x$, since  
$\ax(x_r) + \ax(x_r x_{-1}) = \lambda_r \otimes \lambda_r$,
$\ax(x_r) - \ax(x_r x_{-1}) = 4 X_r$. The images of $98280_x$
under $\tau^{\pm 1}$ cover $98304_x$; hence the images of an axis
span $\Griess$.

So it suffices to show that an axis $v$ has only finitely many images
under the automorphism group. Otherwise there would be an axis $v'$ very
close to $v$  that we can write as $v' = v + \theta w + O(\theta^2)$,
with $w$ a nonzero vector orthogonal to $v$ and $\theta$ small. 
Since $ v * v = 16 v$, we obtain
\begin{align*}
v'  * v' &= 16 \, v + 2\theta\, v * w + O(\theta^2) \, , \\
16 \, v' &= 16 \, v + 16  \theta\, w + O(\theta^2) \, .
\end{align*}
Since $|v * w| < 4\, |w|$ by Lemma~\ref{Lemma:pre:central:2A:axis}, these two equations prohibit
$v'  * v' = 16\, v'$ for sufficiently small $\theta$.  \hfill \proofend

\medskip

Recall that $x_{-1}$  is the (unique) central involution in  $\C$.

\begin{Definition}\rm
A 2B {\em involution} of the Monster is an element conjugate in $\M$ to~$x_{-1}$.
\end{Definition}

From the construction, the trace of a 2B involution on $\Griess$ is $276$;
by Lemma~\ref{Lemma:pre:central:2A:axis}, the trace of a 2A involution on $\Griess$ is $4372$. 
Thus, 2A and 2B involutions are not conjugate in $\M$. 
We will see later in Theorem~\ref{theorem:classification} that every involution of the Monster 
is either a 2A or a 2B involution.

\begin{Theorem}[see~\cite{Conway:Construct:Monster}, §~13]
	\label{Theorem:centralizer:z}	
	Let $x_{-1}$ be the central involution in 
	$\C \cong 2^{1+24}_+.\Co_1$. The centralizer of $x_{-1}$ 
	in $\Aut(\Griess)$ and hence in $\M $ is $\C$. \hfill\proofend
\end{Theorem}

\medskip

For computing the order of $\M$ we use two specific
2A involutions $\beta^+ = x_{ij}$ with $i = 2$, $j = 3$, and
$\beta^- = \beta^+ \cdot x_{-1}$.  We have 
$\beta^+$, $\beta^- \in \Qx0$. We abbreviate $\lambda_{\beta^+}$
to $\lambda_\beta$. Then $\lambda_{\beta^-} =  \lambda_{\beta^+}  = \lambda_\beta \in \Lambda / 2 \Lambda$.
We define the axes $v^+$, $v^-$ introduced in Section~\ref{section:Introduction}
by $v^+ = \ax(\beta^+)$, $v^- = \ax(\beta^-)$. Then the group
$\M_{v^+} \cap \M_{v^-}$ fixing both axes, $v^+$ and $v^-$,
also fixes the 2B involution $x_{-1}$. Hence
$\M_{v^+} \cap \M_{v^-}$ is isomorphic to the subgroup of $\C$
of structure $2^{1+23}.\Co_2$ fixing $\beta^+$.
The mmgroup package provides support for computing
with axes $v^+$, $v^-$, since they are used for shortening
a word of generators of~$\M$.
%%%%%%%%%%%%%%%%%%%%%%%%%%%%%%%%%%%%%%%%%%%%%%%%%%%%%%%%%

\subsection{The Monster and the Moonshine Module}

In this section, we revisit the construction of the Monster 
from the more advanced perspective of vertex operator algebras. 
This allows us to provide more 
conceptual proofs of the results from the last section
and sometimes to obtain stronger results.
We assume here that the reader is familiar with vertex operator algebras
and the construction of the Moonshine module.

%%%%%%%%%%%%%%%%%%%%%%%%%%%%%%%%%%%%%%%%%%%%%%%%%%%%%

\paragraph{The Monster representation.}
The Moonshine module $\V$ is a \voa of central charge~$24$ having
as conformal character
$\chi_{\V}=q^{-1}\sum_{n=0}^{\infty}{\rm dim}\, \V_n\, q^n$ the elliptic
modular function $j-744=q^{-1}+ 196884\, q +  21493760\,q^2 +\cdots $.
It was constructed in~\cite{FLM} (see also~\cite{Bo-ur}) as a
${\bf Z}_2$-orbifold of the lattice \voa $V_{\Lambda}$
associated to the Leech lattice $\Lambda$. 
The algebra induced on the degree-two part $\V_2$ can be 
identified with the Griess algebra $\Griess$ (see~\cite{FLM}, Prop.~10.3.6).

The $-1$ automorphism of $\Lambda$ has an up-to-conjugation unique lift $x_{-1}$ (cf.~\cite{DGH-virs}) to $V_{\Lambda}$ whose 
${\bf Z}_2$-orbifold $V_\Lambda^+\oplus V_\Lambda^{T,+}$ is the Moonshine module.
It follows from the construction that the group $\C\cong \Ciso$
acts by automorphisms on $\V$ such that the subspaces $V_\Lambda^+$ and $V_\Lambda^{T,+}$ 
are the $+1$ and $-1$ eigenspaces of the action of the central element $x_{-1}$ of~$\C$.

\medskip

The Leech lattice $\Lambda$ can be constructed from the Golay code, which allows one to construct the 
additional triality automorphism $\tau$
of $\V$. If $\Golay\subset \F^{24}$ is
the Golay code, then the index~$2$ sublattice $\Lambda_0$ of $\Lambda$ consists of the vectors $\frac{1}{\sqrt{2}}v\in \R^{24}$ such that
$v$ has integer coordinates and the sum of the $24$ coordinates of $v$ is even. The fixed-point vertex operator algebra
$V_\Lambda^{++}:=V_{\Lambda_0}^+\subset V_{\Lambda}^+ \subset \V$ has $16$ irreducible modules, and the associated modular tensor category 
$\mathcal{T}(V_\Lambda^{++})$ is isomorphic to the modular tensor category associated to the direct sum $ H\oplus H $ of two hyperbolic planes~\cite{DongModules,DongFusion}.

For the automorphism group of $V_\Lambda^{++}$ we have the short exact sequence
$$ 1 \longrightarrow \Aut(V_\Lambda^{++})_0  \longrightarrow \Aut(V_\Lambda^{++})  \longrightarrow \overline{\Aut}(V_\Lambda^{++})  \longrightarrow 1, $$
where  $\overline{\Aut}(V_\Lambda^{++})\subset O(H\oplus H)$ is the induced permutation group permuting 
the irreducible modules of $V_\Lambda^{++}$.
The $16$ elements of $ H\oplus H $ form three orbits under the orthogonal group $O(H\oplus H)$ of order $96$:
two orbits of isotropic vectors of size $1$ and $9$, and one orbit of size $6$.
The $16$ modules of $V_\Lambda^{++}$ have $4$ different conformal characters: the orbit of size $9$ splits into two orbits of size $3$ and $6$
under $\overline{\Aut}(V_\Lambda^{++})$. More precisely, $\overline{\Aut}(V_\Lambda^{++})$ is a dihedral group of order $12$, and
$\Aut(V_\Lambda^{++})_0\cong 2^{22}.(2^{11}{:}M_{24})$; cf.~\cite{Shimakura}.

\smallskip

There are three orbits of $2$-dimensional isotropic subspaces $A$ in  $H\oplus H$ under
$\overline{\Aut}(V_\Lambda^{++})$: Extending $V_\Lambda^{++}$ by the simple currents in $A$ results in one of the three self-dual
vertex operator algebras $V_{N({A_1^{24}})}$ (the lattice vertex operator algebra associated to the Niemeier lattice with root system $A_1^{24}$), $V_\Lambda$, or $\V$.
The setwise stabilizer of $A$ in $\overline{\Aut}(V_\Lambda^{++})$ for the case $\V$ is an $S_3$.
It can be shown that $\Aut(\V)$ has a subgroup  $ \N0\cong 2^2.2^{11+22}.(M_{24} \times S_3)$, where
the $2^2$ is the dual group $\hat A$ of $A$ and the $S_3$ is the induced action on $A$; cf.~\cite{FLM}.

\medskip

One checks that the definitions of $\C$ and $\N0$ on $\V$ agree with the definitions in Section~\ref{Monster:Griess} when restricted to the Griess algebra $\Griess=\V_2$. Let $\M'=\langle \C, \N0 \rangle\subseteq \Aut(\V)$. It was shown in~\cite{FLM} that $\M'=\M$. The result also follows from the more general Theorem~\ref{theorem:aut-griess-moonshine} proven below.
This allows us to define the Monster $\M$ also as the group of automorphisms of $\V$ generated by $\C$ and $\N0$.

%%%%%%%%%%%%%%%%%%%%%%%%%%%%%%%%%%%%%%%%%%%%%%%%

\paragraph{Involutions and Ising vectors.}  We study involutions in $\M$ and $\Aut(\V)$ from the viewpoint of vertex operator algebras.
Recall that to an automorphism of order~$n$ of a self-dual \voa one can assign its {\em type} $h \pmod{n}$, where $h|n$, cf.~\cite{EMS20}.

\begin{Lemma}\label{VOA:Lemma:characters}	
Let $t$ be an involution in $\Aut(\V)$ of type~$0$. 
Then the vector-valued character ${\rm Ch}_{(\V)^t}$ of the fixed-point \voa $(\V)^t$ is given either by
$${\rm Ch}_{(\V)^t}=\renewcommand{\arraystretch}{1.2}
\left(\begin{array}{c}	
q^{-1} + 100628\, q + O(q^2) \\
  96256\, q + O(q^2) \\
 96256\, q + O(q^2) \\
q^{-1/2} + 4372\, q^{1/2}+ O(q^{3/2}) 
\end{array}\right)
\quad or \qquad 
\renewcommand{\arraystretch}{1.2}
\left(\begin{array}{c}	
q^{-1} + 98580\, q + O(q^2) \\
24+  98304\, q + O(q^2) \\
 98304\, q + O(q^2) \\
4096\, q^{1/2} + O(q^{3/2}) 
\end{array}\right)_.
$$
\end{Lemma}

\fatline{Proof:}
The general formula for the vector-valued character of the fixed-point \voa $V^t$ of a type~$0$ involution on a self-dual \voa $V$ of central
charge $c=24$ is given in~\cite{HoMo2}, Section~3.5. Since $\V_1=(\V_1)^t=0$, it follows that the three parameters are $a=24$, $b=24-\dim(\widetilde{\V}_1)$ (where $\widetilde{\V}$ is the orbifold \voa of $\V$ with respect to $t$), and $\ell=b/24$. Also, $\ell$ is the dimension of a vector space, i.e., has to be a non-negative integer.
Hence, the only two possibilities are $b=24$ or $b=0$, and the result follows.
\hfill\proofend

\begin{Definition}\rm
An involution $t$ of $\Aut(\V)$ of type~$0$ is called a 2A\nobreakdash-{\em like involution} respectively a 2B\nobreakdash-{\em like involution} if 
the vector-valued character of $(\V)^t$ is the first respectively the second of the two characters given in Lemma~\ref{VOA:Lemma:characters}.
\end{Definition}
We see from Lemma~\ref{VOA:Lemma:characters} that the trace of a 2A like involution on the Griess algebra~$\Griess$ is
$4372$ and that the trace of a 2B like involution is $276$.
We will see in the last section that all involutions $t$ in $\Aut(\V)$ are of type~$0$. For the 2A and 2B involutions of the Monster, this 
can be checked explicitly using a decomposition of $\V$ with respect to a Virasoro frame.
%Therefore one has that the 2A involutions of the Monster are 2A like involutions
%and 2B involutions of the Monster are 2B like involutions of $\Aut\V$.

\medskip

\begin{Definition}\label{Definition:Ising}\rm
An {\em Ising vector} of a vertex operator algebra $V$ is a vector $x\in V_2$ which generates 
a simple Virasoro vertex operator subalgebra of central charge $c=1/2$.
\end{Definition}
A vector $x\in V_2$ generates a Virasoro vertex subalgebra if 
$$x_1x=2\cdot x, \qquad x_2x=0, \quad \hbox{and} \qquad x_3x = \frac{c}{2}\cdot {\bf 1}.$$
The third equation is equivalent to $\langle x,x \rangle =\frac{c}{2}$,
where $\langle \,.\,,\,.\,\rangle $ is the natural invariant bilinear form on $V$.
In the case of the Moonshine module $\V$, one has $x_2x\in \V_1=0$ so that the second 
equation is automatically satisfied. 
Also, because $\V$ has a positive definite invariant
bilinear form, the Virasoro vertex operator subalgebra is automatically simple, i.e.,
isomorphic to the minimal model $L(\sfrac{1}{2},0)$. Thus a vector $x$ in $\Griess$ is an Ising
vector if $x* x =  2\, x$ and $\langle x,x \rangle =\frac{1}{4}$.

\begin{Theorem}\label{VOA:Theorem:Ising-2A}
There exists a natural bijection between 2A like involutions in $\Aut(\V)$ and Ising vectors in $\V$.
\end{Theorem}

\fatline{Proof:} This is a special case of the correspondence between involutions of type~$0$ in the automorphism group of self-dual vertex operator algebras of
central charge $8d$ and self-dual vertex operator superalgebras of charge less than or equal to $8d$ as established in~\cite{HoMo2}.

In the case of a 2A like involution $t$ in $\Aut(\V)$, the vertex operator superalgebra extension $V$ of the fixed-point subalgebra $(\V)^t$ has as character
the sum of the first and last component of the vector-valued modular function: 
$$\chi_V=q^{-1}+q^{-1/2}+ 4372\, q^{1/2} +  100628\, q + O(q^{3/2}).$$
In particular, one has $\ell = \dim V_{1/2}=1$. It follows that the vertex operator subalgebra $\langle V_{1/2} \rangle$ of $V$ is isomorphic to the
``single Fermion'' $L(\sfrac{1}{2},0)\oplus L(\sfrac{1}{2},\sfrac{1}{2})$, and thus one canonically has 
the vertex operator subalgebra $L(\sfrac{1}{2},0)\subset (\V)^t \subset \V$.
The Virasoro element of $L(\sfrac{1}{2},0)$ is the associated Ising vector $x$ in $\V$.

Conversely, given an Ising vector $x$ in $\V$, it generates a subalgebra $L(\sfrac{1}{2},0)\subset \V$. 
One has a natural decomposition (cf.~\cite{Ho-dr})
$$\V=VB^\natural_{(0)}\otimes L(\sfrac{1}{2},0) \oplus VB^\natural_{(1)}\otimes L(\sfrac{1}{2},\sfrac{1}{2})
\oplus VB^\natural_{(2)}\otimes L(\sfrac{1}{2},\sfrac{1}{16}). $$
Then $t$ is defined as the involution acting by $+1$ on the first two terms of this sum and by $-1$ on the third. From the fusion rules of 
$L(\sfrac{1}{2},0)$, it follows that $t$ is an automorphism of $\V$ called the Miamoto involution associated to $x$, and 
one sees that $t$ is of type~$0$.
From the known characters of the $L(\sfrac{1}{2},h)$ for $h=0$, $\sfrac{1}{2}$, and $\sfrac{1}{16}$ and the character of $VB^\natural= VB^\natural_{(0)}\oplus VB^\natural_{(1)}$
as an extremal self-dual vertex operator superalgebra of central charge $c=23\frac{1}{2}$, one sees that $t$ has to be 2A like.

One also checks that the provided maps between 2A like involutions and Ising vectors are inverse to each other.
\hfill\proofend

\smallskip
The VOA $L(\sfrac{1}{2},0)$ has three isomorphism classes of irreducible modules: $L(\sfrac{1}{2},0)$, 
$L(\sfrac{1}{2},\sfrac{1}{2})$, and $L(\sfrac{1}{2},\sfrac{1}{16})$, of conformal weight
$0$, $\sfrac{1}{2}$, and $\sfrac{1}{16}$, respectively. It follows that the endomorphism $x_1={\rm ad}\,x $ of $V_2$ has eigenvalues 
$0$, $\sfrac{1}{2}$, $\sfrac{1}{16}$, and~$2$.

\medskip
The connection to 2A and 2B involutions in the Monster as defined in the previous subsection is explained by the following two results.

\begin{Lemma}\label{VOA:Lemma:2A}	
Let $t$ be a 2A involution of the Monster. Then $t$ is a 2A like involution of $\Aut(\V)$, and the axis associated to $t$ is,
up to a scaling factor, an Ising vector of $\V$.
\end{Lemma}
\fatline{Proof:} An axis $v$ in $\Griess$ was normalized by $v*v=16\,v$ and $(v,v)=8$. 
Thus $x=v/8$ satisfies $x*x=2\,x$ and $(x,x)=1/8$. The scalar product $(\,.\,,\,.\,)$ on $\Griess$ is normalized such that for the identity element $1_\Griess$ one has
$(1_\Griess,1_\Griess)=3/2$, i.e., $(\omega,\omega)=6$ for the Virasoro element $\omega=2\,1_\Griess$ of $\V$, whereas $\langle \omega,\omega \rangle =12$. Hence 
$\langle x,x\rangle=1/4$ and $x$ is an Ising vector. 

By Lemma~\ref{Lemma:pre:central:2A:axis}, the involution associated to $v$ negates the eigenspace of $\ad v$ for the eigenvalue $1/2$.
Thus $\ad x$ negates the eigenspace for the eigenvalue $\sfrac{1}{16}$ in agreement with the definition 
for the associated 2A like involution  
given in the proof of Theorem~\ref{VOA:Theorem:Ising-2A}. \hfill \proofend

\smallskip

We will see in the last section that every 2A like involution of $\Aut(V)$ is actually a 2A involution of the Monster, or,
equivalently, that every Ising vector in $\V$ is, up to scaling, an axis in $\Griess$.

\medskip

\begin{Theorem}\label{VOA:Theorem:2B}	
Let $t$ be a 2B like involution of $\Aut(\V)$. Then $t$ is conjugated in $\Aut(\V)$ to a 2B involution of the Monster,
and the centralizer of $t$ in $\Aut(\V)$ is conjugated to $\C$.
\end{Theorem}

\fatline{Proof:} The fixed-point \voa $(\V)^t$ of an involution of type~$0$ has as its modular tensor category the tensor category associated to a 
hyperbolic plane as the quadratic space. This allows one to construct the orbifold \voa $\widetilde{\V}$ by extending $(\V)^t$ with
the $(\V)^t$ module of conformal weight $0\pmod{1}$ not contained in $\V$. From Lemma~\ref{VOA:Lemma:characters}, it follows that $\widetilde{\V}$ has
the character $q^{-1} + 24 + 196884\, q + O(q^2)$.
It is well known that the Leech lattice \voa $V_\Lambda$ is, up to isomorphism, the unique self-dual \voa $V$ of central charge $c=24$ with
$\dim V_1=24$. Since a lift of the $-1$ involution on $\Lambda$ to $V_\Lambda$ is, up to conjugacy, the unique involution of $V_\Lambda$
having $\V$ as its orbifold (cf.~\cite{HoMo2}, Table~9, entry \#968), one has  $\widetilde{\V}\cong V_\Lambda $ and hence $(\V)^t\cong V_\Lambda^+$. 
It follows that $t$ is conjugated to $z$ in $\Aut(\V)$.

Using~\cite{Shimakura}, we have that the fixed-point \voa $V_\Lambda^+$ has the automorphism group $2^{24}.\Co_1\cong \C/\langle z\rangle $.
Since the centralizer of $z$ in $\Aut(\V)$ has to respect the decomposition of $\V$ into the $z$ eigenspaces $V_\Lambda^+$ and $V_\Lambda^{T,+}$, 
it follows from Schur's Lemma that the centralizer of $z$ is contained in $\mathbb{C}^*.2^{24}.\Co_1 $. Using the fusion rules for $V_\Lambda^+$, it follows 
that the centralizer is equal to $\C\cong \Ciso$.
\hfill \proofend

\smallskip

We note that Theorem~\ref{VOA:Theorem:2B} provides a \voa proof of Lemma~\ref{Theorem:centralizer:z}.

%%%%%%%%%%%%%%%%%%%%%%%%%%%%%%%%%%%%%%%%%%%%%%%%%%%%%

\subsection{Identification of $\Aut(\Griess)$ and  $\Aut(\V)$}

We will show in Section~\ref{section:order:monster} that the Monster is the full automorphism group of the Moonshine module and the Griess algebra. For this, we will need the following general result.
\begin{Theorem}\label{theorem:aut-griess-moonshine}
The map $\rho: \Aut(\V)\longrightarrow \Aut(\Griess)$, which restricts a \voa automorphism of $\V$ to an algebra automorphism of its degree-two part $\V_2=\Griess$, is an isomorphism. 
\end{Theorem}
Here, an automorphism of $\Griess$ is understood to also respect the invariant bilinear form on~$\Griess$.

For the proof, one could try to proceed as in~\cite{FLM} and relate the algebra structure on $\Griess$ (together with the invariant bilinear form) and
the vertex operator algebra structure on $\V$ by utilizing the {\em affine Griess algebra} $\hat\Griess$.
%Since $\hat\Griess$ is defined in terms of $\V$, the $\hat\Griess$ module $\V$ is
%actually an $\Aut(\V)$-equivariant module. Thus the action of $\Aut(\V)$ on $\Griess$ is also faithful and 
%one has $\Aut(\V)\subset \Aut(\Griess)$.
Instead, we make use of the axes in $\Griess$.
A system of $48$ pairwise-orthogonal axes was first found in~\cite{MeNe}.
This motivated the notion of a {\em Virasoro frame} of a vertex operator algebra in~\cite{DGH-virs}.

\smallskip

A  Virasoro frame of a \voa $W$ of central charge $c=n/2$ is a set $\{\omega_1,\,\ldots,\,\omega_n\}\subset W_2$ 
of mutually orthogonal Ising vectors with  $\omega_1+\cdots+\omega_n=\omega$, the Virasoro element of $W$. Such a Virasoro frame generates a vertex operator subalgebra $T_n\cong L(\sfrac{1}{2},0)^{\otimes n}$ inside a sufficiently regular $W$, where $L(\sfrac{1}{2},0)$ is the $c=1/2$ Virasoro minimal model vertex operator algebra.
It was shown in~\cite{DGH-virs} that such a subalgebra $T_{n}$ in the \voa $W$ defines two linear binary codes ${\cal C}$ 
and ${\cal D}$ of length $n$ which are orthogonal to each other. 

\medskip

\fatline{Proof of injectivity of $\rho$:} 
A Virasoro frame $\{\omega_1,\,\ldots,\,\omega_{48}\}\subset \V$ was found in~\cite{DGH-virs}
such that ${\cal C}$ is the lexicographic code of length $48$ and minimal weight~$4$.
The dimension of ${\cal C}$ is $41$ and the code ${\cal D}={\cal C}^\perp$ has dimension~$7$ and minimal weight $16$. 
Furthermore, ${\cal C}$ is generated by its weight~$4$ codewords.

\smallskip

Let $\psi$ be an automorphism of $\V$ in the kernel of~$\rho$. Since $\psi$ fixes $\Griess$, it fixes the Virasoro frame and hence the 
vertex operator algebra $T_{48}$ it generates.
 
\smallskip

The simple vertex operator subalgebra $M_{\cal C} \cong \bigoplus_{c \in{\cal C}} M(c)$ of $\V$ is a 
simple current extension of $T_{48}$ by the simple currents $M(c)$, $c\in {\cal C}$, and the 
fusion rules among the $M(c)$ are those for the group ${\cal C}$. 
A codeword $c$ of weight~$4$ in ${\cal C}$  corresponds to a $T_{48}$-module $M(c)$ with a lowest-weight vector $u$ in $\V_2=\Griess$.
Thus $\psi$ fixes such $u$ and the $T_{48}$-module $M(c)$.
Since ${\cal C}$ is generated by its weight~$4$ codewords, $\psi$ acts trivially on all of $M_{\cal C}$.

Similarly, $\bigoplus_{I\in {\cal D}} M^I \cong \V$ is a simple current extension. 
The $M^I$  are irreducible $M_{\cal C}$-modules, and the fusion rules among the $M^I$,  $I\in {\cal D}$, are those for the group ${\cal D}$. Furthermore, $ {\cal D}$ is generated by codewords $I\in  {\cal D}$ with the property that the $M_{\cal C}$-module~$M^I$ contains a lowest-weight vector $u$ in $\V_2$, cf.~\cite{DGL-uniqueness}.  This implies that $\psi$ acts trivially on those~$M^I$ and thus on all of $\V$.
\hfill \proofend

\medskip

For the proof of the surjectivity of $\rho$, we use the following result from~\cite{DGL-uniqueness}, which was also proven with the help of the above Virasoro frame.
\begin{Theorem}[\!\!{\cite{DGL-uniqueness}, Thm.~1}]
Let $V$ be a $C_2$-cofinite self-dual extremal \voa of central charge~$24$ such that the degree-two part $V_2$ with the induced algebra structure is isomorphic to the Griess algebra $\Griess$. Then $V$ is isomorphic to $\V$.
\end{Theorem}

An analysis of the proof shows that the following more precise result can be deduced:
\begin{Proposition}\label{proposition:aut-voa-extension}
Let $V$ be a $C_2$-cofinite self-dual extremal \voa of central charge~$24$. Let $\varphi: \Griess\longrightarrow V_2$ be an algebra isomorphism.
Then $\varphi$ can be extended to a \voa isomorphism $\hat \varphi: \V\longrightarrow V$.
\end{Proposition}

\fatline{Proof:}
Let $\{\varphi(\omega_1),\,\ldots,\,\varphi(\omega_{48})\}$ be the image of the Virasoro frame $\{\omega_1,\,\ldots,\,\omega_{48}\}$ of 
$\Griess$ under $\varphi$. This is a Virasoro frame in $V_2$ because $\varphi$ is an algebra isomorphism, and thus preserves the product relations defining a frame.
It generates a vertex operator subalgebra $T_{48}'\cong L(\sfrac{1}{2},0)^{\otimes {48}}\subset V$, see~\cite{DGL-uniqueness}. This allows us to extend $\varphi$ from $ \{\omega_1,\,\ldots,\,\omega_{48}\}\subset \V_2 $ to a \voa isomorphism $\hat \varphi :T_{48} \longrightarrow T_{48}'$.

The vertex operator subalgebra $T'_{48}$ of $V$ defines the two associated binary codes ${\cal C}'$ and ${\cal D}'$. 
The codewords $c$ of weight~$4$ in  ${\cal C}$  correspond to $T_{48}$-modules $M(c)$ with a lowest-weight vector $u$ in 
$\V_2=\Griess$. Since $\varphi$ is an algebra isomorphism, it preserves the correspondence between weight~$4$ codewords and their associated lowest-weight vectors in the degree-two space. This establishes a canonical identification between the weight~$4$ codewords of 
${\cal C}$ and ${\cal C}'$. As both codes are generated by their weight~$4$ codewords,  $\varphi$ induces
a canonical identification between ${\cal C}$ and $ {\cal C}'$.
The two vertex operator subalgebras 
$$M_{\cal C} \cong \bigoplus_{c \in{\cal C}} M(c) \qquad \hbox{and} \qquad  M_{\cal C'} \cong \bigoplus_{c \in{\cal C'}} M(c) $$
of $\V$ and $V$, respectively, are therefore isomorphic since they are both simple current extensions of $L(\sfrac{1}{2},0)^{\otimes {48}}$ by corresponding simple currents. 
This allows us to extend the given algebra isomorphism $\varphi:\bigoplus_{c \in{\cal C}} M(c)\cap \Griess\longrightarrow \bigoplus_{c \in{\cal C'}} M(c)\cap V_2$ to a \voa isomorphism $\hat\varphi: M_{\cal C}\longrightarrow M_{\cal C'}$.

It is shown in~\cite{DGL-uniqueness}, Thm.~6.10, that all irreducible $M_{\cal C}$-modules
are simple currents.  It follows that one has simple current extensions
$$\V \cong \bigoplus_{I\in {\cal D}} M^I  \qquad \hbox{and} \qquad V \cong \bigoplus_{I\in {\cal D'}} M^I,$$
where $M^I$ is an irreducible $M_{\cal C}$-module corresponding to a codeword $I\in {\cal D}$ and the fusion rules among the $M^I$ are the ones for the group ${\cal D}$. 
Analogous properties hold for the $M_{\cal C'}$-modules.
Furthermore, $ {\cal D}$ is generated by codewords $I\in  {\cal D}$ with the property that $M^I$ contains an $M_{\cal C}$-lowest-weight vector $u$ in $\V_2$. 
Since $\varphi$ is an algebra isomorphism, these vectors $u$ are mapped to $M_{\cal C'}$-lowest-weight vectors $\varphi(u)$ in  $V_2$. 
%Both extensions are isomorphic vertex operator algebras since they are both equivalent simple current extensions of $M_{\cal C}$ and $M_{\cal C'}$ by corresponding simple currents. 
The code $ {\cal D'}={\cal C'}^\perp$ is, like ${\cal D}$, generated by codewords $I'$ for which the corresponding $M_{\cal C'}$-modules 
$M^{I'}$ contain a lowest-weight vector in $V_2$. The isomorphism $\varphi$ 
ensures that the set of such generating modules for $\V$ corresponds exactly to that for $V$.
This allows us to extend the \voa isomorphism $\hat\varphi : M_{\cal C}\longrightarrow M_{\cal C'}$ to a \voa isomorphism  $\hat\varphi :\V\longrightarrow V$.
\hfill \proofend

\medskip

Taking $V=\V$ in Proposition~\ref{proposition:aut-voa-extension} shows that any automorphism of $\Griess$ 
can be extended to an automorphism of $\V$.
Therefore, the map $\rho$ is surjective. Together with the already proven injectivity, this completes the proof of 
Theorem~\ref{theorem:aut-griess-moonshine}.

%%%%%%%%%%%%%%%%%%%%%%%%%%%%%%%%%%%%%%%%%%%%%%%%%%%%%%%%%%%%%%%%%%%%%%%%%%%%%%%%%%%%%%%

\section{Counting the axes in the Griess algebra}
\label{section:2A:axes:monster}

In this section, we compute the number of axes in the Griess algebra $\Griess$ by decomposing them into twelve orbits under the group $\C$.
The result is shown in Table~\ref{table:orbits:G_x0_axes}. In particular, as the number of all axes in $\Griess$ is the sum of the sizes of all orbits, we will obtain:
\begin{Proposition}\label{prop:numberaxes}
There are $97\,239\,461\,142\,009\,186\,000$ axes in the Griess algebra $\Griess$.
\end{Proposition}
The number of axes in each orbit was first determined by Norton, cf.~\cite{Norton98}, Table~2, where the relevant 
information is given without proof. The main objective of this section is to recompute
the information in that table using the mmgroup package~\cite{mmgroup2020}
without assuming any further properties of the Monster or its character table.

\begin{table}[ht]
	\centering
	\small % Reduces font size for large tables
	\renewcommand{\arraystretch}{1.08} % Increases row height
\begin{tabular}{|c|r|r|r|l|l|}
\hline
$\C$-  &  \multicolumn{2}{c|}{No.\ suborbits} 
 & \multicolumn{1}{c|}{$\C$-orbit size} &
  \multicolumn{2}{c|}{ Stabilizer $C$ of an axis in the orbit}  \\	
\cline{2-3} \cline{5-6}
orbit & $N_{x0}$ & $N_{xyz}$ &  &
  $ \! C  \cap  \Qx0 \!$ &  $C / (C \cap \Qx0)$   \\
\hline  
2A & 3 & 3 & 196560 &  $2^{1+23}$  & $\mbox{Co}_2$  \\
2B & 5 & 6 & 11935123200 &  $2^{1+16}$  & $2^{1+8}.O_8^+(2)$  \\
4A & 9 & 11 & 1630347264000 &  $2^{12}$  & $2^{11}.M_{23}$  \\
4B & 16 & 21 & 1466587938816000 &  $2^{7}$  & $2^{1+8}.S_6(2)$  \\
4C & 14 & 17 & 6599645724672000 &  $2^{5}$  & $2^{1+8+6}.A_8$  \\
6A & 10 & 14 & 1896194506752000 &  $2^{2}$  & $U_6(2).2$  \\
6C & 23 & 35 & 438020931059712000 &  $2^{2}$  & $2^{1+8}.(3 \times U_4(2)).2$  \\
8B & 37 & 64 & 8601138282627072000 &  $1$  & $2^{1+10}.M_{11}$  \\
6F & 8 & 12 & 1501786049347584000 &  $1$  & $2^{1+8}.A_9$  \\
10A & 23 & 36 & 786389785840189440 &  $2$  & $\mbox{HS}.2$  \\
10B & 38 & 68 & 37845008443559116800 &  $1$  & $2^{1+8}.(A_5 \times A_5).2$  \\
12C & 65 & 118 & 48057153579122688000 &  $1$  & $2.S_6(2)$  \\
\hline
Total & 251 & 405 & 97239461142009186000 \\
\cline{1-4}
\end{tabular}
\caption{Sizes of $\C$-orbits on axes and stabilizers in $\C$ of 
	an axis in the orbit.}  
\
\label{table:orbits:G_x0_axes}	
\end{table}

\smallskip

We also use the MAGMA algebra system \cite{bosma1997magma} 
for obtaining the structure of certain groups. 
%%%%%%%%%%%%%%%%%%%%%%%%%%%%%%%%%%%%%%%%%%%%%%%%%%%%%%%%%%%%%%%%%%%%%%%%%%%%%%%%%%%%

\subsection{The strategy for counting the axes} % for counting the 2A axes}
\label{subsection:basic:axes}

\paragraph{Identifying the $\C$-orbits.}

We can find at least twelve orbits.
\begin{Lemma}\label{Lemma:traces}
There exist at least twelve orbits on the axes under $\C$, distinguished by
the conjugacy class of $t x_{-1}$ in $\M$ where $t$ is the
2A involution associated to an axis.
\end{Lemma}
\fatline{Proof:} Starting with the axis $v^+$, we compute samples of
axes by repeatedly transforming an axis with
a random element of $\C$ and with a power of the triality element $\tau$.
The first transformation does not change the $\C$-orbit of the
axis, but the second transformation may do so. In principle, we could
compute traces of powers of  $t x_{-1}$ in $\Griess$ with the mmgroup
package for distinguishing between conjugacy classes in $\M$.
But for reasons of speed, we use a watermarking technique explained below for
distinguishing between orbits of axes. This way we obtain
representatives of twelve different orbits of axes. The Python script
{\tt orbit\_classes.py} computes orders and traces in  $\Griess$
of the elements  $t x_{-1}$ (and also of some powers of $t x_{-1}$),
where $t$ runs through the involutions corresponding to the axes found. 
That way, twelve different conjugacy classes of elements $t x_{-1}$
have been found. 
\hfill\quad  \proofend

\smallskip

\noindent{\bf Remark:} In \cite{Norton98} and \cite{Seysen22},
the $\C$-orbit of a 2A involution $t$
is labeled by the conjugacy class of $t x_{-1}$ in the Monster. The character
information computed by the script {\tt orbit\_classes.py} in the accompanying 
code is sufficient to match our representatives of the orbits on the axes with the
orbits described in \cite{Norton98}. Here the character table of the Monster
in the ATLAS~\cite{Atlas} is only used for adjusting our notation
for the orbits to the notation used in \cite{Norton98}.

\medskip

We choose representatives for the axes for each of the twelve 
cases as in Lemma~\ref{Lemma:traces}.
To confirm that a given axis belongs to one of the corresponding twelve
$\C$-orbits, we proceed as follows using a watermarking technique to
identify the likely orbit.

The subspace $300_x$ of $\Griess$ can be
identified with the symmetric $24 \times 24$ matrices. 
An element of $\C$ acts as a rational orthogonal
transformation on the symmetric $24 \times 24$ matrices.
Thus the eigenvalues and dimensions of the eigenspaces of 
such a matrix depend therefore only on the $\C$-orbit.
We project the given axis to $300_x$, obtaining a symmetric $24 \times 24$ matrix.
Table~1 in \cite{Seysen22} shows that the algebraic information obtained 
in that way is sufficient to distinguish between the putative twelve  $\C$-orbits.
A slightly more sophisticated disambiguation of the $\C$-orbits of the
axes based on their entries in $\Griess$ (modulo~$15$)
is given in~\cite{Seysen22}, Section~8.4.

Having identified the likely orbit of an axis, we want to show that
the axis is actually in that orbit. Therefore, 
we compute an element of $\C$ that maps the axis to the corresponding
representative. For the details of this computation, we refer
to Appendix~\ref{appendix:map:orbit:Gx0}.

\medskip

\paragraph{The triality transition matrix.}

Our objective is to compute the transition matrix~$M$ for the 
operation of the triality element $\tau$ on the $\C$-orbits,
as shown in Table~\ref{table:action:G_x0_axes}. The
rows and columns of $M$ are labeled by the names of the $\C$-orbits on the axes. 
An entry in row $i$, column $j$ indicates the number of axes in orbit $j$ 
that are mapped to orbit $i$ by the triality element $\tau$, up to a factor depending on
the column $j$ only and chosen such that the total of each column is
$16584750$. Thus the matrix $M/16584750$ is column-stochastic,
and one checks that it is regular.

From the Perron–Frobenius theorem, we conclude that
the column vector containing the sizes of the $\C$-orbits is the unique eigenvector
of $M$ for the eigenvalue~$16584750$. For
computing the sizes of the $\C$-orbits on the axes from the
matrix~$M$, it, therefore, suffices to know the size of one of these orbits. 
\begin{Lemma}
The orbit labeled '2A' and containing the axis $v^+$ has size $196560$.
\end{Lemma}
\fatline{Proof:}
The 2A~involution $\beta^+$ corresponding to the axis $v^+$ in the orbit
labeled '2A' is contained in $\Qx0$. The element $\lambda_\beta$ of
$\Lambda/2\Lambda$ is short and has precisely the two preimages
$\beta^+$ and $\beta^-$ in $\Qx0$. Since there are $98280$ short
vectors in  $\Lambda/2\Lambda$, the size of the  $\C$-orbit
on the axes labeled '2A' is $2 \cdot 98280$. 
\hfill \proofend

\smallskip

So given $M$, we may compute the sizes of all $\C$-orbits on the axes, as
shown in Table~\ref{table:orbits:G_x0_axes}.

\begin{table}[ht]
	\centering
	\small % Reduces font size for large tables
	\renewcommand{\arraystretch}{1.08} % Increases row height
	
	% First Array
	\begin{tabular}{|c|r|r|r|r|r|r|}
		\hline
		& 2A & 2B & 4A & 4B & 4C & 6A \\
		\hline
		2A & 93150 & 135 & 1 & . & . & . \\
		2B & 8197200 & 64935 & 3542 & 63 & 15 & . \\
		4A & 8294400 & 483840 & 35927 & 2232 & 240 & 891 \\
		4B & . & 7741440 & 2007808 & 137367 & 33600 & 74844 \\
		4C & . & 8294400 & 971520 & 151200 & 54255 & . \\
		6A & . & . & 1036288 & 96768 & . & 133731 \\
		6C & . & . & 7254016 & 2225664 & 860160 & 2020788 \\
		8B & . & . & 5275648 & 3870720 & 430080 & 8382528 \\
		6F & . & . & . & . & 522240 & . \\
		10A & . & . & . & 2064384 & 491520 & 953856 \\
		10B & . & . & . & 1548288 & 7311360 & . \\
		12C & . & . & . & 6488064 & 6881280 & 5018112 \\
		\hline
	\end{tabular}
	\bigskip
	
	% Second Array
	\begin{tabular}{|c|r|r|r|r|r|r|}
		\hline
		& 6C & 8B & 6F & 10A & 10B & 12C \\
		\hline
		2A & . & . & . & . & . & . \\
		2B & . & . & . & . & . & . \\
		4A & 27 & 1 & . & . & . & . \\
		4B &  7452 & 660 & . & 3850 & 60 & 198 \\
		4C & 12960 & 330 & 2295 & 4125 & 1275 & 945 \\
		6A & 8748 & 1848 & . & 2300 & . & 198 \\
		6C & 321003 & 134024 & . & 284900 & 59520 & 72450 \\
		8B & 2631744 & 2217999 & 483840 & 3010000 & 1375680 & 1399104 \\
		6F & . & 84480 & 486135 & . & 271200 & 274320 \\
		10A & 511488 & 275200 & . & 411775 & 131520 & 106992 \\
		10B & 5142528 & 6052992 & 6834240 & 6329400 & 6435735 & 6543936 \\
		12C & 7948800 & 7817216 & 8778240 & 6538400 & 8309760 & 8186607 \\
		\hline
	\end{tabular}
\captionsetup{justification=centering}
\caption{Action of the triality element $\tau$ on the $\C$-orbits in the axes.
\newline \small
 The entry in row $i$, column $j$
is the number of $\Nxyz$-orbits in the $\C$-orbit $j$ that
\newline are mapped to $\C$-orbit	$i$. The total of each column is $|\C / \Nxyz| = 16584750$.
}  
\label{table:action:G_x0_axes}	
\end{table}

\medskip

To compute the transition matrix $M$, we utilize that the triality element 
$\tau$ normalizes a large subgroup of $\C$. 
The group $\Nx0 = \C \cap \N0$ has a subgroup $\Nxyz$ of index~$2$ with
structure $2^{2+11+22}.M_{24}$. The triality element $\tau$ normalizes
$\Nxyz$. Thus, for any $g \in \C$ and axis $v$, the $\C$-orbit
of $v g \tau$ depends on the left coset $g \Nxyz$ only. The index
of $\Nxyz$ in~$\C$ is $16584750$. Thus a left transversal 
$\mathcal{T}$ of $\Nxyz$ in  $\C$ can be effectively computed.
For computing the transition matrix, it suffices to transform a
representative of each $\C$-orbit with all elements
$\{g \cdot \tau \mid g \in  \mathcal{T}\}$, and 
to determine the $\C$-orbit of the transformed axes.

Such a computation is lengthy, but doable with a modern PC, and
has been done in an early stage of the mmgroup project using the 
identification of an orbit via watermarking. We have
$\M = \langle \C, \tau \rangle.$ Therefore, if every orbit obtained
by watermarking in this way is one of the twelve $\C$-orbits given by
Lemma~\ref{Lemma:traces}, then this calculation also shows 
that we have found all $\C$-orbits.
% since the triality element $\tau$ maps the  $\C$ orbits into each other.

\medskip

It thus remains to verify that the orbit of an axis identified by watermarking
is indeed in that $\C$-orbit. Here, in principle, the method discussed in
Appendix~\ref{appendix:map:orbit:Gx0}
can be used. However, we have to transform almost 200 million axes to the representative 
of their $\C$-orbit. Finding a suitable transformation takes much longer than identifying
an orbit via watermarking, so this approach is presently not feasible.
In the next section, we will discuss how to reduce the
required amount of computation significantly.

\medskip

Once the data for matrix $M$ have been computed, their correctness
can be verified by a much simpler computation. For details, we
refer to Appendix~\ref{appendix:certificate}.

%%%%%%%%%%%%%%%%%%%%%%%%%%%%%%%%%%%%%%%%%%%%%%%%%%%%%%%%%%%%%%%%%%%%%%%%%%%%%%%%%%%%%%%%%%%%%%%%

\subsection{Speeding up the counting of the axes} % for counting the 2A axes}
\label{subsection:speedup:axes}

In this subsection, we explain how to speed up the counting of the axes so
that this can be done in practice. Our primary objective is
to compute Table~\ref{table:action:G_x0_axes}.	
Let $v_i$ be the chosen representative of the $\C$-orbit $i$ on the axes, 
and let $C_i$ be the stabilizer of $v_i$ in $\C$. 
As discussed in the previous subsection,
we want to compute the axes $v_i g \tau$, 
where~$g$ runs through a left transversal $\mathcal{T}$
of $\Nxyz$ in $\C$. 

Since $C_i$ centralizes $v_i$, it
suffices to consider a transversal of the double coset space
$C_i \backslash \C / \Nxyz$, provided that we know the
number of right cosets of $C_i$ in $\C$ in each double
coset in $C_i \backslash \C / \Nxyz$. So we need a fast
method to find the transversal of the double coset space
$C_i \backslash \C / \Nxyz$ and to count the right cosets of
$C_i$ in each of these double cosets.
Therefore, it would be helpful to find a permutation representation
of  $\C$ acting on the cosets in  $\C/ \Nxyz$.
The group  $\C$ acts on $\Lambda / 2 \Lambda$ via the
natural action of its factor group $\Co_1$. The
centralizer of $\lambda_\Omega$ (with $\lambda_\Omega$ as in
Section~\ref{Monster:Griess}) in $\C$ is $\Nx0$. Since
$\lambda_\Omega$ is a type 4 vector, and $\Co_1$ is transitive on
the type 4 vectors, there is a one-to-one correspondence between
the left cosets  $\C / \Nx0$ and the type 4 vectors
in $\Lambda / 2\Lambda$. Since $\Nxyz$ has index~$2$ in $\Nx0$,
this is almost what we need.

There is an $a \in  \Nx0\setminus \Nxyz$ with 
$a^{-1} \tau a = \tau^{-1}$; we may take, e.g., 
$a = x_\delta$ with $\delta = \{0\} \in \mathcal{C^*}$ in the notation of 
\cite{Seysen20}. So instead of inspecting the representatives 
$v_i g \tau$ and $v_i g a \tau$ of two different cosets in
$ \C / \Nxyz$, we may inspect $v_i g \tau$ and $v_i g \tau^{-1}$
instead, where $g$ runs through a transversal of
$C_i \backslash \C / \Nx0$.

\smallskip

We first model the stabilizer $C_i$ in $\C$ of an axis $v_i$,
or at least a large subgroup of $C_i$.
For this, we create a set of generators
of $C_i$ as elements of the group $\C$. The easiest way to create
a generator of $C_i$ is to multiply $v_i$ with a random element $g_1$
of $\C$, and to find an element $g_2 \in \C$ that reduces the axis
$v_i g_1$ to $v_i$, i.e., $(v_i g_1) g_2 = v_i$. Then $g_1 g_2$
centralizes $v_i$. We may obtain a suitable element $g_2$ by 
the method discussed in Appendix~\ref{appendix:map:orbit:Gx0}.
That way, we have computed 10 random generators of the stabilizer
$C_i$. There is a very small probability
that these generators generate a proper subgroup of $C_i$.
For computing Table~\ref{table:action:G_x0_axes}, this may cause
some more overhead, but not any error. For simplicity, we also
write  $C_i$ for the group generated by these generators.

We model the action of the subgroup  $C_i$ of $\C$ on
$\Lambda / 2\Lambda$ by using the class {\tt Orbit\_Lin2} in the
mmgroup package. This class implements the action of a general
group $C_i$ (given by generators) on the vector space $\mathbb{F}_2^n$,
$n \leq 24$, or, more specifically, on $\Lambda / 2\Lambda$. 
The key ingredient of the class {\tt Orbit\_Lin2} is a fast 
algorithm for finding orbits on $\Lambda / 2\Lambda$ under the group $C_i$.
That class contains a member function for computing a list of
representatives of
the orbits under $C_i$ on $\Lambda / 2\Lambda$ and the sizes of these
orbits. Since $C_i \subset \C$, all elements of such an orbit
have the same type in $\Lambda / 2\Lambda$, and we are only
interested in orbits on vectors of type~$4$. 

For any $\lambda \in \Lambda / 2\Lambda$ of type 4, we may compute
a $g_\lambda \in \C$ that maps $\lambda$ to  $\lambda_\Omega$
with mmgroup.
Let $L_i$ be a set of representatives of all orbits of type 4
vectors under $C_i$ on $\Lambda / 2\Lambda$. Then the set
$ \{g_\lambda \mid  \lambda \in L_i\}$ is a transversal
of the double coset space  $C_i \backslash \C / \Nx0$.
We compute the  $\C$-orbits on the axes
$v_i g_\lambda \tau^{\pm 1}$, where $\lambda$ runs through
$L_i$, as discussed  in Section~\ref{subsection:basic:axes}. In
the course of this computation, we also verify that all these axes
are in one of the twelve orbits given by Lemma~\ref{Lemma:traces}.
This proves:

\begin{Lemma}
There are precisely twelve orbits on the axes under $\C$.
\end{Lemma}
	 
We compute the column of Table~\ref{table:action:G_x0_axes}
corresponding to the orbit $v_i \C$ as follows.	For each 
$\lambda \in L_i$ and $\varepsilon = \pm1$, we add the size of
the orbit of $\lambda$ under $C_i$ on $\Lambda / 2\Lambda$ to
the row in the table corresponding to the orbit of the axis
$v_i g_\lambda \tau^{\varepsilon}$.

Having computed Table~\ref{table:action:G_x0_axes}, we may compute
column~4 of Table~\ref{table:orbits:G_x0_axes}, as discussed in
Section~\ref{subsection:basic:axes}; in particular, Proposition~\ref{prop:numberaxes}
is proven.

%%%%%%%%%%%%%%%%%%%%%%%%%%%%%%%%%%%%%%%
\subsection{Decomposition of the $\C$ orbits into $\Nxyz$ orbits}
\label{subsection:reduce:no:orbits}

Let the axis $v_i$ be a representative of the orbit  $v_i \C$ as above.
Once the size of that orbit is known (as shown in column~4
of Table~\ref{table:orbits:G_x0_axes}), we may compute 
the decomposition of the orbit $v_i \C$ into $\Nx0$-orbits
and also the structures of the stabilizers in $\Nx0$ of 
representatives of these orbits.

In the previous subsection, we have computed a set of random
elements of the stabilizer $C_i$ of $v_i$ in $\C$. There is
a very small probability that this set generates a proper
subgroup of $C_i$. For computing the $\Nx0$-orbits
on the axes in the orbit $v_i \C$, we need a set generating
the whole group $C_i$.
In Appendix~\ref{appendix:compute:centralizers:axes},
we will show how to check that a computed set of elements
of $C_i$ actually generates $C_i$. In the
remainder of this section, we assume that a generating set
of the whole group  $C_i$ is available.

If $C_i$ is the stabilizer of an axis $v_i$ in $\C$, then the number
of $\Nx0$-orbits on the axes contained in the orbit  $v_i\C$ is equal
to the number of double cosets in $C_i \backslash \C / \Nx0$. 
As we have seen in the previous subsection, the left cosets
$\C / \Nx0$ correspond to the type~4 vectors in $\Lambda / 2\Lambda$.
Thus, the double cosets in $C_i \backslash \C / \Nx0$ correspond to
the orbits of the type~4 vectors in $\Lambda / 2\Lambda$ under $C_i$.
These orbits can easily be enumerated with the technique used
in the previous subsection, so we may compute column~2 in
Table~\ref{table:orbits:G_x0_axes}. Clearly, the total of that
column is equal to the number of orbits on axes under $\Nx0$,
and hence also to the number of double cosets 
in $\Mv \backslash \M / \Nx0$. The stabilizer $\M_{v^+}$ is of structure $2.B$.
Thus, we have shown:
 
 \begin{Proposition}
 	There are precisely 251 double cosets in $2.B \backslash \M / \Nx0$.
 	The number of double cosets in 	$2.B \backslash \M / \Nx0$
 	contained in a double coset  in 	$2.B \backslash \M / \C$
 	is given by Table~\ref{table:orbits:G_x0_axes}, Column~2.		
 	\label{prop:double:2B:M:Nx0}
 \end{Proposition}
 
 We remark that the number of double cosets in
 $2.B \backslash \M / \Nx0$ (with  $\Nx0$ of structure
 $2^{2+11+22}.(M_{24} \times 2)$)
 has not been computed before.

The last two columns of Table~\ref{table:orbits:G_x0_axes} describe
the structure of the stabilizer in $\C$ of an axis in a given
$\C$-orbit. The computation of this structure is described
in Appendix~\ref{appendix:compute:centralizers:axes}.
We also list the number of double cosets in $2.B \backslash \M / \Nxyz$ in
column~3 for each of the twelve $\C$-orbits. The computation of these
numbers is described in Appendix~\ref{N0:orbits}.

\bigskip

We have verified all the entries of Table~\ref{table:orbits:G_x0_axes}.
We also determined the  $\N0$-orbits on the axes, and
the result is shown in Appendix~\ref{N0:orbits}.

%%%%%%%%%%%%%%%%%%%%%%%%%%%%%%%%%%%%%%%%%%%%%%%%%%%%%%%%%%%%%%%%%%%%%%%%%%%%%%%%%%%%%%%

\section{Counting the feasible axes}
\label{section:compute:centralizers:axes:2B}

In this section we compute the number of feasible axes.
Recall that $\beta^+$ and $\beta^-$ are the 2A involutions in $\M$
corresponding to the specific axes $v^+$ and $v^-$ as 
defined in Section~\ref{Monster:Griess}.
We define the set of  {\em feasible axes}, denoted $X^-$, as
the orbit of the axis $v^-$ under the stabilizer $\Mv$ of axis~$v^+$,  i.e.
$$X^-=\{v^-h \mid h \in \M_{v^+}\}.$$
We write $H$ for the stabilizer $\Mv \cap \C$ of  $v^+$ in $\C$, as in
\cite{Seysen22}. The group $H$ is also the centralizer of
$\beta^+$ in $\C$, and of structure $2^{1+23}.\mbox{Co}_2$.
Since the triality element $\tau$ centralizes $\beta^+$, 
it must also fix its corresponding axis $v^+$ and 
we have $\tau \in \M_{v^+}$.
Our strategy for counting the feasible axes is analogous to the method
for counting the axes in Section~\ref{section:2A:axes:monster}:
we first decompose $X^-$ into orbits under the action of $H$
and then compute a transition matrix describing how $\tau$ acts on these orbits.
We will just summarize this briefly. We will show:

\begin{Proposition}\label{prop:numberbabyaxes}
There are $11\,707\,448\,673\,375$ feasible axes.
\end{Proposition}

%Since $v^+$ is orthogonal to $v^-$ the whole set $X^-$ consists of axes orthogonal to the axis $v^+$.
%However, we do not show or use the fact that $X^-$ is the set of {\em all} axes orthogonal to the axis $v^+$.
%Let $\lambda_\beta$ The 2A~involution $\beta^+$ corresponding to the axis $v^+$ in the orbit
%labeled with '2A' is contained in $\Qx0$. The image 
%$\lambda_\beta$ of the involution $\beta^+$ in $\Lambda/2\Lambda$

We use a method similar to that in  Section~\ref{section:2A:axes:monster}
for disambiguating orbits on the feasible axes under $H$.
We adopt the naming conventions for these orbits from~\cite{Seysen22}, Section~9.

%We have $\beta^+$, $\beta^- \in \Qx0$; and we write 
%$\lambda_\beta$ for the common image of $\beta^+$ and $\beta^-$
%in $\Lambda / 2 \Lambda$ under the natural homomorphism
%$\Qx0 \rightarrow   \Lambda / 2 \Lambda$. 
\begin{Lemma}\label{Lemma:baby:traces}
	There exist at least ten orbits on the feasible axes under $H$. 
	These axes are distinguished
	by the conjugacy class of $t x_{-1}$ in $\M$ where $t$ is the
	2A involution associated to an axis, and by the norm
	of the projection of the axis onto the subspace of $300_x$
	spanned by $\lambda_\beta \otimes \lambda_\beta$.
\end{Lemma}

\fatline{Proof:} We generate a large sample of feasible axes by applying random elements of $H$
and powers of the triality element $\tau$ to the initial axis $v^-$.
A similar watermarking strategy as in the proof of Lemma~\ref{Lemma:traces}
allows to disambiguate ten different $H$-orbits on the feasible axes.
Those can be distinguished by the two invariants described in the statement of the lemma.
\hfill \proofend

\medskip

\noindent{\bf Remark:} 
M\"uller \cite{mueller_2008} has enumerated the ten orbits on the
feasible axes under the action of
$H = \Mv \cap \C$ on~$X^-$ and computed their sizes in
\cite{mueller_2008}, Table~3. However, that paper uses properties of
the Baby Monster and its character table.
We will confirm these results with mmgroup and provide some additional
information for each orbit in our Table~\ref{table:orbits:H_axes}.

\medskip

We choose representatives for the feasible axes for each of the ten 
potential orbits found in Lemma~\ref{Lemma:baby:traces}. 
To confirm that a given feasible axis belongs to one of these potential
$H$-orbits we use a watermarking technique to identify the likely orbit,
similar to the watermarking used in Section~\ref{subsection:basic:axes}.
(A slightly more sophisticated method for the disambiguation of the $H$-orbits
of the feasible axes is discussed in~\cite{Seysen22}, Section~9.2.)

Having identified the likely orbit of a feasible axis, 
we compute an element of $H$ that maps the axis to the corresponding
representative, as discussed in Appendix~\ref{appendix:map:orbit:H}.

%%%%%%%%%%%%%%%%%%%%%%%%%%%%%%%%%%%%%%%%%%%%%

%\paragraph{The $H$ orbits on $X^-$.}

\begin{table}[ht]
	\centering
	\small % Reduces font size for large tables
	\renewcommand{\arraystretch}{1.08} % Increases row height
	\begin{tabular}{|c|r|r|r|l|l|}
		\hline
		$H$-  &  \multicolumn{2}{c|}{No.\  suborbits} &
		 \multicolumn{1}{c|}{$H$-orbit size} &
		\multicolumn{2}{c|}{Stabilizer $C$ of an axis in the orbit}  \\	
		\cline{2-3} \cline{5-6}
		orbit & $H \cap N_{x0}$ &  $H \cap N_{xyz}$  & &
		$ \! C  \cap  \Qx0 \!$ &  $C / (C \cap \Qx0)$   \\
		\hline
2A1 & 1 & 1 & 1 &  $2^{1+23}$  & $\mbox{Co}_2$  \\
2A0 & 5 & 6 & 93150 &  $2^{1+22}$  & $2^{10}.M_{22}.2$  \\
2B1 & 4 & 5 & 7286400 &  $2^{1+16}$  & $2^{1+8}.S_6(2)$  \\
2B0 & 6 & 8 & 262310400 &  $2^{1+15}$  & $2^{1+4+6}.A_8$  \\
4A1 & 9 & 12 & 4196966400 &  $2^{12}$  & $2^{9+1}.L_3(4).2$  \\
4B1 & 12 & 18 & 470060236800 &  $2^{7}$  & $2^{1+8+5}.S_6$  \\
4C1 & 6 & 9 & 537211699200 &  $2^{5}$  & $2^{1+4+6}.A_8$  \\
6A1 & 3 & 5 & 9646899200 &  $2^{2}$  & $U_6(2).2$  \\
6C1 & 8 & 14 & 6685301145600 &  $2^{2}$  & $2^{1+8}.U_4(2).2$  \\
10A1 & 5 & 9 & 4000762036224 &  $2$  & $\mbox{HS}.2$  \\
\hline
Total & 59 & 87 & 11707448673375 \\
		\cline{1-4}
	\end{tabular}
	\captionsetup{justification=centering}
	\caption{Sizes of $H$-orbits on feasible axes (with $H = \Mv \cap \C$),  
		\newline 
	    and stabilizers in $H$ of an axis in the orbit. 
	}  
	\
	\label{table:orbits:H_axes}	
\end{table}

%Sizes of $H$ orbits on feasible axes and stabilizers of these axes in $H$

\paragraph{The triality transition matrix for $X^-$.}

Our objective is to compute the transition matrix $M'$
for the operation of the triality element $\tau$ on the
$H$-orbits on the feasible axes, as shown in
Table~\ref{table:action:tau:H}. Matrix $M'$ agrees with
the matrix shown in \cite{mueller_2008}, Table~3. Rows and 
columns of $M'$  correspond to the $H$-orbits on the feasible axes;
entries are interpreted as the entries of matrix $M$ in
Table~\ref{table:action:G_x0_axes}. The columns in the
table are normalized so that they sum to	
$|H{:}H \cap \Nxyz| = 93150$. The computation of $M'$
will be discussed below.

\begin{table}[ht]
	\centering
	\small
	%\scriptsize % Reduces font size for large tables
	\renewcommand{\arraystretch}{1.08} % Increases row height	
	\setlength{\tabcolsep}{4pt}
	\begin{tabular}{|c|r|r|r|r|r|r|r|r|r|r|}
		\hline
		&     2A1 &     2A0 &     2B1 &     2B0 &     4A1 &     4B1 &     4C1 &     6A1 &     6C1 &    10A1 \\
		\hline
		2A1 &       . &       1 &       . &       . &       . &       . &       . &       . &       . &       . \\
		2A0 &   93150 &     925 &      63 &      15 &       1 &       . &       . &       . &       . &       . \\
		2B1 &       . &    4928 &      63 &     120 &      42 &       1 &       . &       . &       . &       . \\
		2B0 &       . &   42240 &    4320 &    1815 &     420 &      30 &      15 &       . &       . &       . \\
		4A1 &       . &   45056 &   24192 &    6720 &    1807 &     272 &     120 &     891 &      27 &       . \\
		4B1 &       . &       . &   64512 &   53760 &   30464 &   10287 &    5040 &   24948 &    3060 &    3850 \\
		4C1 &       . &       . &       . &   30720 &   15360 &    5760 &    3495 &       . &    4320 &    4125 \\
		6A1 &       . &       . &       . &       . &    2048 &     512 &       . &     891 &      36 &     100 \\
		6C1 &       . &       . &       . &       . &   43008 &   43520 &   53760 &   24948 &   53451 &   53900 \\
		10A1&       . &       . &       . &       . &       . &   32768 &   30720 &   41472 &   32256 &   31175 \\
		\hline
	\end{tabular}
	\bigskip
	\captionsetup{justification=centering}
	\caption{Action of the triality element $\tau$ on 
		the $H$-orbits in the feasible axes. 
		\newline \small
		The entry in row $i$, column $j$
		is the number of $(H \cap \Nxyz)$-orbits in $H$-orbit
		\newline $j$  that are mapped to $H$-orbit $i$. 
		The total of each column is $|H{:}H \cap \Nxyz| = 93150$.
	}  
	\label{table:action:tau:H}
\end{table}

Using the same argument based on Markov chains as in
Section~\ref{subsection:basic:axes} we may compute a vector
containing the sizes of the $H$-orbits on
the feasible axes (up to a scalar multiple) from
matrix $M'$. We have:

\begin{Lemma}
   The orbit labeled with '2A1' and containing the axis $v^-$ has size $1$.
\end{Lemma}
\fatline{Proof:}
The group $H$ fixes the two involutions
$\beta^+$ and $x_{-1}$, and hence also their product $\beta^-$
and the axis $v^- = \ax(\beta^-)$. 
Axis $v^-$ is in orbit 2A1; so orbit 2A1 has size~$1$. 
\hfill\proofend

\smallskip

Thus we may compute the sizes of all $H$-orbits on the feasible axes,
as shown in Table~\ref{table:orbits:H_axes}, column~4.

For computing the transition matrix $M'$  we use the
fact that the triality element $\tau$ normalizes the
group $H \cap \Nxyz$ with structure $2^{2+11+20}.M_{22}.2$.
Thus for any $g \in H$ and feasible axis $v$ the $H$-orbit of
$v g \tau$ depends only on the left coset $g(H \cap\Nxyz)$. 
The index $|H{:}H \cap \Nxyz|$ is equal to $93150$; hence a
left transversal $\mathcal{T}'$ of  $H / (H \cap \Nxyz)$ 
in $H$ can be effectively computed. For computing the transition
matrix $M'$ it suffices to transform a representative of each
$H$-orbit with all elements $\{g \cdot \tau \mid g \in \mathcal{T}' \}$
and to determine the $H$-orbits of the transformed axes.

In Appendix~\ref{appendix:map:orbit:H} we will discuss how to
effectively determine the $H$-orbit of a feasible axis. We may
speed up the computation of $M'$ using a similar technique as
described in Section~\ref{subsection:speedup:axes}. We also verify
that the images of all feasible axes obtained in the course
of this computation are in one of the ten orbits given by
Lemma~\ref{Lemma:baby:traces}. This proves:

\begin{Lemma}
	There are precisely ten orbits on the feasible axes under $H$.
\end{Lemma}

Having computed matrix $M'$ in Table~\ref{table:action:tau:H},
we may compute column~4 of Table~\ref{table:orbits:H_axes}, as
discussed above; in particular, Proposition \ref{prop:numberbabyaxes} is proven.

\medskip

For obtaining the last two columns of Table~\ref{table:orbits:H_axes}
we compute the structure of the stabilizer $C_i$ of a feasible axis
$v_i$ in $H$, where $v_i$ runs through the representatives of
the $H$-orbits on the feasible axes. We may compute random
elements of $C_i$ as follows. We multiply $v_i$ by a random
element $g_1$ of $H$. Then we use the technique described in
Appendix~\ref{appendix:map:orbit:H} to compute a $g_2 \in H$ that
maps $v_i g_1$ to $v_i$; i.e.~$g_1 g_2$ centralizes $v_i$.
That way we have computed 10 random elements of the stabilizer
$C_i$. These elements generate the whole group $C_i$ with a very
high probability. The method discussed  in
Appendix~\ref{appendix:compute:centralizers:axes} can be used
to check if a set of elements of $C_i$ actually generates $C_i$.
%So we may assume that a generating system for $C_i$ is available.

So we have a generating system for the stabilizer
in $H$ of a feasible axis, for each $H$-orbit on the feasible axes.
The structures of these stabilizers have been computed with MAGMA and
mmgroup, similar to the corresponding computation in
Appendix~\ref{appendix:compute:centralizers:axes} for the $\C$-orbits. These structures
are displayed in the last two columns of Table~\ref{table:orbits:H_axes}.

%!!!!!!!!!!!!!!!!!!!!!!!!!!!!!!!!!!!!!!!!!!!!!!!!!!!!!!!

It remains to compute columns~2 and~3 of
Table~\ref{table:orbits:H_axes}. As in \cite{Seysen22}, we call
a vector $\lambda \in \Lambda / 2 \Lambda$ {\em feasible} if
$\lambda$ is of type~2 and $\lambda + \lambda_\beta$ is
of type~4. Then $H \cap N_{x0}$ is the stabilizer of the feasible
vector $\lambda_\beta + \lambda_\Omega$, see \cite{Seysen22},
Section~9.1.
Hence there is a one-to-one correspondence between
the left cosets  $H / (H \cap \Nx0)$ and the feasible type~2
vectors in $\Lambda / 2\Lambda$.

Let $C_i$ be the stabilizer of a feasible axis $v_i$ in $H$. Since
the left cosets  $H / (H \cap \Nx0)$ correspond to  the feasible
type~2 vectors in $\Lambda / 2\Lambda$, the double cosets in
$C_i \backslash H /  (H \cap \Nx0)$ correspond to the orbits in
the feasible type~2 vectors in $\Lambda / 2\Lambda$ under $C_i$.
These orbits can easily be enumerated with the technique described
in Section~\ref{subsection:reduce:no:orbits}; so we may compute
column~2 in Table~\ref{table:orbits:H_axes}. Clearly, the total
of that column is equal to the number of orbits on feasible axes under
$H\cap \Nx0$, and hence also to the number of double cosets 
in  $H \backslash \Mv / (H \cap \Nx0)$.  The stabilizer $\Mv$ is of structure
$2.B$. Thus we have shown:

\begin{Proposition}
	There are precisely $59$ double cosets in $H \backslash 2.B / (H \cap \Nx0)$.
	The number of double cosets in \hbox{$H \backslash 2.B / (H \cap \Nx0)$}
	contained in a double coset in \hbox{$H \backslash 2.B / H$}
	is given by Table~\ref{table:orbits:H_axes}, Column~2.
	\label{prop:double:H:2B:Nx0}
\end{Proposition}

It is well known that the group $H$ of structure $2^{1+23}.\mbox{Co}_2$  is
a maximal subgroup of $2.B$, and that the group $H \cap \Nx0$ of structure 
$2^{2+11+20}.(M_{22}:2 \times 2)$ is maximal in $H$, see, e.g., \cite{Atlas}.
We also list the number of double cosets in
$H \backslash 2.B / (H \cap \Nxyz)$ in
column~3 for each of the ten $H$-orbits. The computation of these
numbers is described in Appendix~\ref{N0capH:orbits}.

\bigskip

We have  verified all the entries of Table~\ref{table:orbits:H_axes}.
We also determined the $(2.B\cap \N0)$-orbits  on the feasible axes and
the result is shown in Appendix~\ref{N0capH:orbits}.

%%%%%%%%%%%%%%%%%%%%%%%%%%%%%%%%%%%%%%%%%%%%%%%%%%%%%%%%%%%%%%%%%%%%%%%%%%%%%%%%%%%%%%%%%%%%%%%%%%%%%%%%%%%%%%%%%%%%%%%%%%%%%%%%%%%%
\section{Properties of the Monster}
\label{section:order:monster}

In this section, we use the enumeration results from Sections~\ref{section:2A:axes:monster} 
and~\ref{section:compute:centralizers:axes:2B}
to obtain several results about the Monster.

\smallskip

\begin{Theorem}\label{theorem:order:monster}
The order of the Monster $\M$ equals
$$  2^{46}\cdot 3^{20}\cdot 5^9\cdot 7^6 \cdot  11^2 \cdot  13^3 \cdot 17 
\cdot 19 \cdot 23\cdot  29 \cdot  31\cdot   41 \cdot  47 \cdot  59 \cdot  71 \qquad\qquad \qquad\qquad\qquad $$
$$\qquad  =\  808,017,424,794,512,875,886,459,904,961,710,757,005,754,368,000,000,000. $$
\end{Theorem}

\fatline{Proof:} In Proposition~\ref{prop:numberaxes}, the size of the set $X^+$ of axes was
determined. In Propositions~\ref{prop:numberbabyaxes}, the size of the set $X^-$ of feasible axes 
is given. Thus the order of the Monster is 
$$\M=|X^+| \cdot |X^-|  \cdot |S|$$ 
where $S$ is the stabilizer in $\M$ of the two axes $v^+$ and $v^{-}$.
The product of the two 2A involutions $\beta^+$ and $\beta^-$ associated with $v^+$ and $v^{-}$ is the central 2B involution
$x_{-1}$ in $\C$ which must be centralized by $S$. By Theorem~\ref{Theorem:centralizer:z} or Theorem~\ref{VOA:Theorem:2B}, it follows
that $S$ is contained in~$\C$. The stabilizer of $v^+$ and $v^{-}$ in $\C$ can be read off from Table~\ref{table:orbits:H_axes}
to be $2^{1+23}.\Co_2$. Since the order $42,305,421,312,000$ of $\Co_2$ is assumed to be known, the result follows. \hfill \proofend

\smallskip

Recall that ${\rm C}_{\M}(t)/\langle t\rangle $ where $t$ is a 2A involution in the Monster is the Baby Monster $B$, the second largest of the sporadic simple groups.
Using $|B|=|X^-|\cdot |S| /2$ we have also shown:

\begin{Corollary}\label{theorem:order:babymonster}
The order of the Baby Monster $B$ equals
$$  2^{41}\cdot 3^{13}\cdot 5^6\cdot 7^2 \cdot  11 \cdot  13 \cdot 17 \cdot 19 \cdot 23\cdot   31\cdot   41 \cdot  47  \qquad\qquad\qquad\qquad $$ $$\qquad\qquad\qquad \qquad\qquad = 4,154,781,481,226,426,191,177,580,544,000,000.$$
\end{Corollary}

%%%%%%%%%%%%%%%

\medskip

Next, we will prove that the Monster is simple. 
Our proof is along the lines of the corresponding proof
in~\cite{Carnahan:2023}. For a group $G$ we write $G^\#$
for the set of non-identity elements of $G$. 
If $G$ is a subgroup and $H$ is a subset of an ambient group $M$,
we write $C_G(H)$ for the centralizer of $H$ in $G$.
\begin{Lemma}
\label{lemma:involution:operating}
If $N$ is a normal subgroup of $\M$ and $z$ is an involution in
$\M$ with $C_N(z) = 1$ then $N$ is abelian of odd order; and
conjugation with $z$ in $N$ is equivalent to inversion in $N$.	
\end{Lemma}

\fatline{Proof:}
When the involution $z$ acts by conjugation on $N$,
it has exactly one fixed point, the identity. Thus $|N|$ is odd. Let
$a$, $b \in N$. The mapping $N \rightarrow N$ given by
$a \mapsto a z a^{-1} z$ is injective (and hence bijective), since 
$a z a^{-1} z = b z b^{-1} z$ implies
$b^{-1} a z a^{-1} b = z$, i.e.~$ a^{-1} b \in C_N(z)$. We have
$ (a z a^{-1} z)^z =  (a z a^{-1} z)^{-1}$. Hence $c^z = c^{-1}$ for all
$c \in N$. Since inversion is an automorphism on $N$,
the group $N$ is abelian. 
\hfill \proofend

\begin{Lemma}
	\label{lemma:normal:Gx0:M}
	A normal subgroup $N$ of $\M$ with $|\C \cap N| = 1$ is trivial.
\end{Lemma}

\fatline{Proof:}
We use the argument from the proof of Lemma~4.8 in~\cite{Carnahan:2023}.
Let $N$ be a normal subgroup of $\M$
with $|\C \cap N| = 1$. Then $|C_N(x_{-1})| = 1$.
Let $E$ be the elementary abelian~$2$ subgroup
$ \{ x_{\pm 1},\, x_{\pm \Omega} \}$ of $\Qx0$,
with $x_{\pm 1}$, $x_{\pm \Omega}$ as in
Section~\ref{Monster:Griess}. Since conjugation
with the triality element $\tau$ cyclically exchanges the
three elements of $E^\#$, we have $|C_N(e)| = 1$ for all
$e \in E^\#$. So by Lemma~\ref{lemma:involution:operating},
$|N|$ is odd, and conjugation in $N$ with any $e \in E^\#$ is
equivalent to inversion in $N$. Thus conjugation with the
product of the three elements of $E^\#$ is also inversion
in $N$. Since that product is equal to 1, inversion is the
identity map in $N$.
Since $|N|$ is odd, this implies $|N| = 1$. \hfill \proofend

\begin{Lemma}
	\label{lemma:series:Gx0}
	The group $G_{x0}$ contains the characteristic series
	\[
	1 \vartriangleleft  \langle x_{-1}\rangle  \vartriangleleft
	\Qx0 \vartriangleleft \C \, .
	\]
	This series contains all normal subgroups of $\C$.
\end{Lemma}

\fatline{Proof:} % New AI-generated proof
We may use the argument in the proof of \cite{Carnahan:2023}, Lemma~4.5.
By construction, $\langle x_{-1} \rangle$ is the center
of $\C$, and $\Qx0$ is a normal subgroup with quotient $\C/\Qx0 \cong \Co_1$.
The series is therefore a normal series. That it contains all normal subgroups
of $\C$ follows from the fact that $\Qx0$ is extraspecial,
 $\Co_1$ is simple, and the natural 24-dimensional representation
$\Qx0/\langle x_{-1} \rangle \cong \Lambda/2\Lambda$ of $\Co_1$
is irreducible over $\mathbb{F}_2$. \hfill \proofend

% This version adds a brief justification for why the series is complete,
% making the argument more self-contained. 

\begin{Theorem}
The Monster is a simple group.
\end{Theorem}

\fatline{Proof:}
We adapt the proof of Theorem~4.12 from \cite{Carnahan:2023} to our
notation. Let $N$ be a non-trivial normal subgroup of $\M$.
By Lemma~\ref{lemma:normal:Gx0:M}, the intersection $N \cap \C$ must be a
non-trivial normal subgroup of $\C$. So by Lemma~\ref{lemma:series:Gx0},
the group $N \cap \C$ must be one of the groups
$ \langle x_{-1}\rangle $, $\Qx0$, or $\C$.
Suppose $N \cap \C = \langle x_{-1} \rangle$. This would imply 
$x_{-\Omega} = (x_{-1})^\tau \in N$. Since
$x_{-\Omega} \in \C \setminus  \langle x_{-1}\rangle$,
this contradicts the assumption $N \cap \C = \langle x_{-1} \rangle$. 

Let $d \in \mathcal{P} \setminus \{\pm1, \pm \Omega\}$, with $ \mathcal{P}$
the {\em Parker loop} as in \cite{Conway:Construct:Monster,Seysen20}.
Then $x_d \in \Qx0$, but $x_d^\tau = y_d \in \C \setminus \Qx0$.
Thus $N \cap \C = \Qx0$ would imply $y_d \in N \cap ( \C \setminus \Qx0)$,
a contradiction.   

Hence $\C \subset N$.
Let $g = x_\delta \in \C$, for any odd element $\delta$ of
the Golay cocode. Then $g^\tau g = \tau$, implying $\tau \in N$.
Hence $N = \M$. 

Note that all calculations in this proof can be done inside the group $\N0$.
\hfill \proofend

\medskip

We now show that the Monster is the full automorphism group of the Griess algebra and the Moonshine module.
This result was first obtained by Tits~\cite{ti-monster,ti-monster2} by reducing it to a group-theoretical characterization of the Monster
given by Smith~\cite{Smith}.
Our method utilizes results of Carnahan~\cite{Carnahan:2023} on the order of $\Aut(\V)$ and Borcherds' proof of the 
Conway-Norton Moonshine conjectures~\cite{Bo-lie}.

\begin{Theorem}
\label{theorem:monster:full}
The Monster is the full automorphism group of the Griess algebra and the Moonshine module, i.e., one has
$\M=\Aut(\Griess)=\Aut(\V)$.
\end{Theorem}

\fatline{Proof:} By Theorem~\ref{theorem:aut-griess-moonshine}, we have $\Aut(\Griess)=\Aut(\V)$. In Carnahan~\cite{Carnahan:2023}, it was
proven that $\Aut(\V)$ equals the order of the Monster as in Theorem~\ref{theorem:order:monster} up to a factor which is a power of $11$; 
see~\cite{Carnahan:2023}, Corollary~5.2 and Theorems~5.3, 5.4, and~5.5. 

Thus we need to show that the Monster has no $11$-Sylow subgroup of order $11^3$ or
larger. If this were the case, there would exist a subgroup $H$ of order~$11^3$ in $\Aut(\Griess)$. By Borcherds~\cite{Bo-lie}, the trace of
any element $g$ of $\Aut(\V)$ on $\V$ is a completely replicable modular function. These functions have been classified in~\cite{ACMS}. 
In particular, there is only one completely replicable function of order $11$ and none of order $11^2$. Thus the exponent of the group $H$ 
must be~$11$. 
In addition, the trace of an order~$11$ element $g$ in $\Aut(\V)$ on $\Griess$ equals~$17$, the coefficient of $q^1$ in the unique completely replicable function of order~$11$. For the dimension of the subspace of $\Griess$ fixed by $H$ we obtain
$$\dim \Griess^H=\frac{1}{|H|} \sum_{h\in H} \tr (h|\Griess)  =
\frac{1}{11^3}\bigl(196884+(11^3-1)\cdot 17\bigr)=\frac{1814}{11}$$
which is not an integer. Thus such a subgroup $H$ cannot exist and the $11$-Sylow subgroup of $\Aut(\V)$ must have order $11^2$.

It follows that $|\Aut(\V)|=|\M|$ and hence $\Aut(\V)=\M$. \hfill \proofend

\medskip

We finally show that the Monster has exactly two conjugacy classes of involutions.
\begin{Theorem}\label{theorem:classification}
The Monster has two conjugacy classes of involutions: the 2A involutions corresponding to the axes and with centralizers isomorphic
to a two-fold cover of the Baby Monster and the 2B involutions with centralizers isomorphic to $\C\cong \Ciso$. 
\end{Theorem}
We provide two proofs, a computational one in Appendix~\ref{appendix:classes:involutions:Monster} using mmgroup
and a theoretical one using vertex operator algebras below.

\smallskip

\fatline{Proof:} It follows from Theorem~\ref{theorem:order:monster} that the index of $\C$ in $\M$ is odd and hence that 
any involution of $\M$ is conjugated to an involution in $\C$.

\smallskip

We first observe that there are no involutions of type $1$ in $\Aut(\V)$. Indeed, by~\cite{Bo-lie,CKU} the graded trace of such an involution 
$t$ equals either the graded trace of a 2A like or a 2B like involution. This determines the character of $(\V)^t$ and hence the whole
vector-valued character of $(\V)^t$. In particular, we see that $(\V)^t$ has no modules of conformal weight $1/4$ or $3/4$. Thus
$t$ must be of type~$0$.

\smallskip

We know from Lemma~\ref{VOA:Lemma:characters} that every involution of $\Aut(\V)$ of type~$0$ is either 2A like or 2B like and 
from Theorem~\ref{VOA:Theorem:2B} that every 2B like involution is actually a 2B involution. We also know from Theorem~\ref{VOA:Theorem:Ising-2A} that every 2A like involution corresponds uniquely to an Ising vector of $\V$. 
Any involution in $\C$ is centralized by the central element $x_{-1}$ in $\C$ and thus the corresponding Ising vector for
a 2A like involution is fixed by $x_{-1}$, i.e., is contained in the fixpoint vertex operator subalgebra $V_{\Lambda}^+$ of~$\V$.
The Ising vectors of $V_\Lambda^+$ have been classified in~\cite{LS-Ising}. They either correspond
to norm~$4$ vectors of $\Lambda$ or to embeddings of $E_8(2)$ (the root lattice $E_8$ rescaled such that the minimal vectors have norm~$4$) into $\Lambda$.
One immediately verifies that those two cases correspond to the first two orbits of axes of $\Griess$ in Table~\ref{table:orbits:G_x0_axes} and thus these Ising vectors are axes in $\Griess$, i.e.,\ the 2A like involution is a 2A involution. 
\hfill \proofend

\smallskip

%Theorem~\ref{theorem:classification} in particular implies that any involution of $\V$ is of type~$0$.

\medskip

The Monster is known to have a total of $194$ conjugacy classes of elements forming $172$ algebraic conjugacy classes~\cite{Atlas}.
Besides two classes of order $27$, the classes of elements $g$ can be distinguished by the traces of $g$ and its powers $g^2$ and $g^3$ 
on $\Griess$.
Explicit representatives have been first determined in~\cite{Barraclough:2005:CCR}.
We provide a list of representatives of these classes for mmgroup at the first author's webpage~\cite{Hoehnclasses}. 
Functions for listing and identifying a class are also available on that webpage.
At the moment we are unable 
to distinguish the two conjugacy classes of elements of order $27$.

%%%%%%%%%%%%%%%%%%%%%%%%%%%%%%%%%%%%%%%%%%%%%%%%%%%%%%%%%%%%%%%%%%%%%%%%%%%%%%
%%%%%%%%%%%%%%%%%%%%%%%%%%%%%%%%%%%%%%%%%%%%%%%%%%%%%%%%%%%%%%%%%%%%%%%%%%%%%%%%%%%%%%%%%%%%%%%%%%%%%%%%%%%%%%%%%%%%%%%%%%%%%%%%%%%%%%%%%%%%%%%%%%

\begin{appendices}

%%%%%%%%%%%%%%%%%%%%%%%%%%%%%%%%%%%%%%%%%%%%%%%%%%%%%%%%%%%%%%%%%%%%%%

\section{Mapping an axis to the representative of its orbit}
\label{appendix:map:orbit}

\subsection{Mapping an axis to the representative of its $\C$  orbit}
\label{appendix:map:orbit:Gx0}

Class {\tt Axis} in the mmgroup package models an axis of a 2A
involution. For documentation we refer to \cite{mmgroup_doc},
Chapter {\em Axes of 2A involutions in the Monster}. Here one of the
most important functions is the member function {\tt reduce\_G\_x0} of
class {\tt Axis}. This
function returns an element of the group $\C$ that maps the axis
to the representative of its $\C$-orbit. This process is called
the {\em reduction} of an axis (in $\C$), and will be discussed
in this section. 

We will not prove that {\em every}
axis can be reduced as described here. It suffices that we are able
to reduce all axes occurring during our practical computations. 
Tests with the mmgroup package have shown that our reduction method
has succeeded on hundreds of thousands 
of random axes. A key step in the implementation of the group
operation in $\M$ in mmgroup is the reduction of an axis in
$\M$. Reduction of an axis in $\M$ is simpler
than in $\C$, because we may transform an axis to a 'simpler' 
$\C$-orbit whenever possible. 

\smallskip

We assume that an axis is given as a vector in the representation
$\Griess$ (with entries taken modulo~$15$).
Our first step is to determine the $\C$-orbit of
the axis. Here we use the method in \cite{Seysen22}, Section~8.4
for a quick disambiguation of the twelve orbits. (At this stage
we are not yet sure that there are no more than twelve orbits; but any
attempt to reduce an axis will certainly fail if the axis is not
in one of these orbits.) 

The next step is called the {\em beautifying} of an axis. In this step
we try to transform the symmetric $24 \times 24$ matrix corresponding
to the part $300_x$ of the axis into a nice form
by using a transformation in $\C$. The same operation is
also applied to a candidate for the representative of a
$\C$-orbit of an axis, in order to obtain a more beautiful
representative. Criteria for a beautiful symmetric matrix are:
\begin{itemize}[itemsep=0pt]
	\item[-]
	Few off-diagonal nonzero entries, so that eigenspaces
	are visible or at least easy to compute.
	\item[-]
	The signs of the nonzero off-diagonal elements should
	follow an orderly pattern. 
\end{itemize}
This beautification is implemented in the function
{\tt beautify\_axis} in file  {\tt beautify\_axis.py} in the
mmgroup package. The parts $300_x$ of the representatives of the
$\C$-orbits on the axes can be displayed with the following
Python script:
\begin{verbatim}
    from mmgroup.axes import Axis
    for name, representative in Axis.representatives().items():
        print("Name of axis orbit:", name)
        print("Part 300_x of representative of orbit:")
        representative.display_sym()
\end{verbatim}
In all cases, after beautifying an axis, the part $300_x$ of the axis is equal
to the corresponding part of one of those representatives. Afterwards, 
it suffices in all but one case to transform the axis with an element
of $\Qx0$, in order to obtain the representative of the
$\C$-orbit of the axis. In the case of the orbit '6F', we have to work in an
extension of structure $2^{1+24+11}$ of $\Qx0$ instead.  

\medskip

In all cases, the reduction of the axis was achieved.

%%%%%%%%%%%%%%%%%%%%%%%%%%%%%%%%%%%
\subsection{Mapping a feasible axis to the representative of its $H$ orbit}
\label{appendix:map:orbit:H}

The reduction of a feasible axis under the group $H= \M_{v^+} \cap \C $ follows a procedure entirely analogous to the one described for $\C$-orbits in Appendix~\ref{appendix:map:orbit:Gx0}. Member function {\tt reduce\_G\_x0} of class
{\tt BabyAxis} in the mmgroup package uses the same principles of orbit disambiguation, beautification, and final reduction to map a feasible axis to the pre-computed representative of its $H$-orbit. For details we refer again to the documentation in \cite{mmgroup_doc}, Chapter {\em Axes of 2A involutions in the Monster}.
 
\medskip
 
As in Appendix~\ref{appendix:map:orbit:Gx0}, we will not prove that
{\em every} feasible axis can be reduced. It suffices that
we are able to reduce all axes occurring during the calculation.
In all cases, the reduction of the axis was achieved.

%Here we also abort the program with an error if this is not possible. 

%%%%%%%%%%%%%%%%%%%%%%%%%%%%%%%%%%%%%%%%%%%%%%%%%%%%%%%%%%%%%%%%%%%%%%%%%%%%%%%%%%%%%%%%%%%%%%%%%%%%

\section{Computing stabilizers of axes in $\C$-orbits}
\label{appendix:compute:centralizers:axes}

A random algorithm is called a {\em Monte Carlo} algorithm if it
may compute a wrong result with a (usually) small probability, 
and a {\em Las Vegas}
algorithm if it either succeeds or indicates failure.
In Section~\ref{subsection:speedup:axes} et seq.\ 
we have computed a generating set
of the stabilizer $C_i$ in $\C$ of a given axis $v_i$  using
a Monte Carlo algorithm. Note that this algorithm computes a
generating set of a group $C'_i$ that might be a proper
subgroup of $C_i$ with a very small probability.
In this section we will convert that algorithm to a 
Las Vegas algorithm. In order to detect $C'_i \subsetneq C_i$, 
it suffices to compute the order $|C'_i|$ of the group $C'_i$.
Then we may multiply $|C'_i|$ by the (known)
size of the $\C$-orbit of $v_i$. If this product is equal
to the order of $\C$, then we have $C'_i = C_i$. Here it
suffices to compute a lower bound for $|C'_i|$
that is sharp with high probability.

\smallskip

While computing the order of a subgroup is a standard feature in computer algebra systems 
like MAGMA~\cite{bosma1997magma} that hardly needs a detailed discussion, 
the group $\C$ is too large to be handled directly. Our approach therefore combines the strengths 
of mmgroup for efficient computation within $\C$ with MAGMA's ability to analyze smaller factor groups.
The computer algebra system MAGMA usually
represents a finite group as a permutation group or as a matrix
group over a finite field.
It includes a large number of algorithms from computational group
theory, such as those for computing group orders, identifying
subgroups, and more; see~\cite{holt2005handbook} for a comprehensive
overview.
The group $\C$ is so large that we have not found a satisfactory
way to compute with that group in MAGMA. Computing in the group
$\C$ in the mmgroup package is easy; but that package contains
only a limited set of algorithms from computational group
theory. It turns out that this is sufficient for computing a
generating set for a subgroup $C'_i$ of the stabilizer 
$C_i$ of an axis, and also a lower bound for
the order of that subgroup. 
In all our tests (over a dozen runs), the computed lower bound for $|C'_i|$ was equal to the expected order of $C_i$
for all cases.

\medskip

As a subgroup of $\C \cong 2^{1+24}.\Co_1$, the group
$C_i$ has a normal subgroup $C_i \cap \Qx0$ with factor group
$\overline{C_i} = C_i / (C_i \cap \Qx0)$. We use the mmgroup package
for representing the generators of $\overline{C_i}$ (with
$\overline{C_i} \subset \mbox{Co}_1 \subset \mathop{\rm SL}_{24}(2)$)
as $24 \times 24$ matrices over $\mathbb{F}_2$.  Then we use
MAGMA to compute the structure of the group 
$\overline{C_i}$ generated by these matrices.
The structure of the group  $\overline{C_i}$ obtained in this way is
displayed in  Column~5 of Table~\ref{table:orbits:G_x0_axes}.
The structure of the subgroup $C_i \cap \Qx0$  can easily be
computed with mmgroup; it is displayed in
Column~4 of Table~\ref{table:orbits:G_x0_axes}.

\medskip

In the remainder of this subsection we will discuss the computation
of a lower bound for $|C'_i|$, with $C'_i \subset \C$ given by
a set of generators. In Section~\ref{subsection:speedup:axes}
we have computed the action of such a group   $C'_i$  on the
type 4 vectors in $\Lambda / 2\Lambda$ by using the class
{\tt Orbit\_Lin2} in the
mmgroup package. This way we obtain a permutation representation
of  $C'_i$ with kernel $C'_i \cap \Qx0$, i.e. a faithful
representation of $\overline{C'_i}$.
The class {\tt Orbit\_Lin2} also supports the computation of
a Schreier vector for the action of the group $C'_i$ on the
type 4 vectors in $\Lambda / 2\Lambda$ as discussed in
\cite{holt2005handbook}, Section~4.1.1. Using a Schreier vector,
we may compute an element of $C'_i$ that maps a vector
in $\Lambda / 2\Lambda$ to any vector in the same $C'_i$-orbit.
So, in principle, a Schreier vector may be used for computing
the stabilizer of any type~4 vector  in $\Lambda / 2\Lambda$.
For simplicity, we assume that this type~4 vector is the vector 
$\lambda_{\Omega}$, as defined in Section~\ref{Monster:Griess}.
Let  $C'_{i,1}$ be the stabilizer
of $\lambda_{\Omega}$ in $C'_i$. Then  $|C'_i|$ is the product
of  $|C'_{i,1}|$ and the size of the orbit of $\lambda_{\Omega}$
under the action of $C'_i$. Since we need just a lower bound
for $|C'_i|$,  it suffices to compute a set of random elements of
$C'_i$ that generates $C'_{i,1}$ with high probability. Since
$C'_{i,1}$ fixes  $\lambda_{\Omega}$, it is a subgroup of the group
$\Nx0\cong 2^{1+24+11}.M_{24}$.  

Let $\rho_{24}$ be the representation of $M_{24}$ as a subgroup
of $\mathop{\SL}_{11}(2)$  given by the action of  $M_{24}$ on
the even part of the Golay cocode $\mathcal{C^*}$. Since the
factor group $C'_i / (C'_i \cap \Qx0)$  is a subgroup of $M_{24}$,
the restriction of representation  $\rho_{24}$ to that factor
group is also a (possibly non-faithful) representation of 
$C'_{i,1}$. We also write $\rho_{24}$ for that representation
of $C'_{i,1}$; and the mmgroup package supports computation
in $\rho_{24}$. Considering $\rho_{24}$ as a permutation
representation of  $C'_{i,1}$ in $\mathbb{F}_2^{11}$, we may
compute the stabilizer of a vector 
$v \in \mathbb{F}_2^{11}$ that is not fixed by $C'_{i,1}$,
using the same method as above. This way we obtain
a subgroup  $C'_{i,2} \subsetneq  C'_{i,1}$;
and the index of  $C'_{i,2}$ in $C'_{i,1}$ is the size of the
orbit of $v$ under the action of $C'_{i,1}$. Again, it
suffices to find a set of elements of $C'_{i,2}$ that 
generates  $C'_{i,2}$ with high probability. Repeating this
process we obtain a chain of subgroups	
\[
   C'_i \supset C'_{i,1} \supset \ldots  \supset C'_{i,k},
\]
such that $C'_{i,k}$ fixes all vectors in $\rho_{24}$.
Thus $C'_{i,k}$ is contained in the normal subgroup of
structure  $2^{1+24+11}$ of $\Nx0$. The mmgroup package
supports computation in the subgroup $2^{1+24+11}$ of $\Nx0$, 
so that the order of $C'_{i,k}$ can be computed.

\smallskip

On the authors' PCs the computation of
Tables~1 and~2 (excluding the computation of the structure of the
stabilizers of the axes with MAGMA) takes about a minute.

%%%%%%%%%%%%%%%%%%%%%%%%%%%%%%%%%%%%%%%%%%%%%%%%%%%%%%%%%%%%%%%%%%%%%%%%%%%%%%

\section{Orbits on axes under $\N0$ and $\N0 \cap 2.B$}

%%%%%%%%%%%%%%%%%%%%%%%%%%%%%%%%%%%%%%

\subsection{Orbits on the axes in $\Griess$ under $\N0$}
\label{N0:orbits}

In Table~\ref{tab:N0:orbits} we list the $123$ orbits on the 2A axes under
the action of the maximal subgroup $\N0$ of $\M$. Recall that $\N0$ has
structure $2^{2+11+22}.M_{24}.S_3$, and that the normal subgroup of $\N0$
of structure $2^{2+11+22}.M_{24}$ is denoted by $\Nxyz$.

\medskip

For each $N_0$-orbit on the axes, we display the structure of the
stabilizer in $\N0$ of an axis in that orbit.
This stabilizer is presented in the form $C = E.G.S$, where
$E = C \cap 2^{2+11+22}$, $E.G = C \cap \Nxyz$. More specifically, we have:
\begin{itemize}
\item $E$ is the intersection of $C$ with the normal 2-subgroup $2^{2+11+22}$ of $\N0$.
It is displayed in Column~2 of the table. 
The structure of $E$ is given as $2^{\lambda+\mu+\nu}$, meaning that its intersections with the characteristic 
subgroups $2^2$ and $2^{2+11}$ of $2^{2+11+22}$ have orders $2^\lambda$ and $2^{\lambda+\mu}$, respectively, 
and the total order of $E$ is $2^{\lambda+\mu+\nu}$.
Here we omit leading zero terms in the exponent describing $E$.  
 
\item $G$ is the projection of $C \cap N_{xyz}$ to the factor group
$N_{xyz}/2^{2+11+22} \cong M_{24}$. Its structure and order
are displayed in Columns~3 and~4. Here Column~3 has been computed
with the GAP computer algebra package~\cite{GAP4}.

\item $S$ is the projection of $C$ to the factor group
$\N0 / \Nxyz \cong S_3$. It is displayed in Column~5. 
\end{itemize}

\smallskip

Since $\N0 / \Nxyz \cong S_3$, an $\N0$-orbit of an axis decomposes
into $|S_3| : |S| = 6 : |S|$ different $\Nxyz$-orbits of the same size,
with $S$ as above.
Since $\Nxyz < \Nx0 < \N0$ and $|\N0:\Nx0| = 3$,  each $\N0$-orbit
decomposes into one, two, or three $\Nx0$-orbits. 
The number of $\Nx0$-orbits in a given $\N0$-orbit is determined by
the action of the group $S$ on the three cosets of $\Nx0$ in $\N0$.
\begin{itemize}
\item If $|S| = 1$, $S$ fixes all three cosets, resulting in three distinct
$\Nx0$-orbits, each of them containing two $\Nxyz$-orbits.
\item If $|S| = 2$, $S$ fixes one coset and permutes the other two, resulting
in two $\Nx0$-orbits. One of them contains one, and the other contains two 
$\Nxyz$-orbits.
\item If $|S| \geq 3$, $S$ acts transitively on the three cosets,
resulting in a single $\Nx0$-orbit which may contain one or two  $\Nxyz$-orbits.
\end{itemize}
The final column, `$\C$-orbits' of Table~\ref{tab:N0:orbits}, describes how
these constituent $\Nx0$-orbits are contained in the twelve $\C$-orbits
from Table~\ref{table:orbits:G_x0_axes}. For example:
\begin{itemize}
\item An entry `$4B$, $6A$, $6C$' indicates that the $\N0$-orbit decomposes
into three $\Nx0$-orbits lying in those $\C$-orbits.
\item An entry `$2A$, $2B^2\mspace{0.8mu}$' indicates that the $\N0$-orbit
decomposes into two $\Nx0$-orbits. The smaller one lies in $\C$-orbit $2A$, 
while the larger one lies in $\C$-orbit $2B$.
\item An entry `$2A^3\mspace{0.8mu}$' indicates that the $\N0$-orbit
contains a single $\Nx0$-orbit lying in the $\C$-orbit $2A$.
\end{itemize}

\smallskip

\smallskip

\ifdim \the\dimexpr\pagegoal-\pagetotal\relax < 20ex
\pagebreak
\fi

\nopagebreak
\vspace{-5ex}
\nopagebreak
{\small 
	\begin{longtable}{|r|l|l|r|c|l|}
	\caption{Orbits of ${N_0}$ on axes} \label{tab:N0:orbits} \\
		\hline
		\multicolumn{6}{|c|}
		{$\N0$-orbits on the axes, and stabilizers
			$\vphantom{2^{2^2}}C = E.G.S$   of these axes in $\N0$ } \\
		\hline
		\textbf{No} & $E$ & $G$ & $|G|$ & $S$ & $\C$-orbits \\
		\hline
		\endfirsthead
		
		\hline
		\textbf{No} & $E$ & $G$ & $|G|$ & $S$ & $\C$-orbits \\
		\hline
		\endhead
		
		\hline
		\endfoot
		
		\hline
		\endlastfoot
		$\vphantom{2^{2^2}}$
		$1$ & $2^{2+11+20}$ & $M_{22} :2$ & $887040$ & $S_3$ & $2A^3$\\
		$2$ & $2^{2+10+16}$ & $2^{4} : A_{8}$ & $322560$ & $2$ & $2A,2B^2$\\
		$3$ & $2^{1+11+11}$ & $M_{23}$ & $10200960$ & $2$ & $2A,4A^2$\\
		$4$ & $2^{2+7+14}$ & $2^{1+6} : \mbox{PSL}_{3}(2)$ & $21504$ & $S_3$ & $2B^3$\\
		$5$ & $2^{1+7+12}$ & $2^{6} : A_{5} : S_{3}$ & $23040$ & $2$ & $2B,4A^2$\\
		$6$ & $2^{6+10}$ & $2^{4} : A_{7}$ & $40320$ & $S_3$ & $4A^3$\\
		$7$ & $2^{1+5+11}$ & $2^{4} : S_{6}$ & $11520$ & $2$ & $2B,4B^2$\\
		$8$ & $2^{1+4+11}$ & $2^{1+6} : \mbox{PSL}_{3}(2)$ & $21504$ & $2$ & $2B,4C^2$\\
		$9$ & $2^{6+6}$ & $2^{4} : A_{8}$ & $322560$ & $2$ & $4A,4B^2$\\
		$10$ & $2^{3+10}$ & $2^{6} : \mbox{PSL}_{3}(2)$ & $10752$ & $S_3$ & $4C^3$\\
		$11$ & $2^{4+8}$ & $2^{1+6} : \mbox{PSL}_{3}(2)$ & $21504$ & $2$ & $4A,4C^2$\\
		$12$ & $2^{2+9}$ & $\mbox{PSL}_{3}(4) :2$ & $40320$ & $2$ & $4A,6A^2$\\
		$13$ & $2^{4+10}$ & $2^{4} : \mbox{PSL}_{3}(2)$ & $2688$ & $2$ & $4A,4B^2$\\
		$14$ & $2^{11}$ & $M_{11}$ & $7920$ & $2$ & $4A,8B^2$\\
		$15$ & $2^{2+9}$ & $2^{4} :3: S_{5}$ & $5760$ & $2$ & $4A,6C^2$\\
		$16$ & $2^{2+8}$ & $2^{5} : S_{5}$ & $3840$ & $S_3$ & $4B^3$\\
		$17$ & $2^{3+6}$ & $2^{6} : (S_{3} \times S_{3})$ & $2304$ & $S_3$ & $4B^3$\\
		$18$ & $2^{3+8}$ & $2^{2+6} : S_{3}$ & $1536$ & $2$ & $4B^2,4C$\\
		$19$ & $2^{1+0}$ & $M_{22}$ & $443520$ & $S_3$ & $6A^3$\\
		$20$ & $2^{1+8}$ & $2^{2+6} : S_{3}$ & $1536$ & $S_3$ & $4C^3$\\
		$21$ & $2^{1+0}$ & $M_{22} :2$ & $887040$ & $2$ & $6A,10A^2$\\
		$22$ & $2^{1+7}$ & $2^{4} : S_{6}$ & $11520$ & $1$ & $4B,6A,6C$\\
		$23$ & $2^{1+6}$ & $2^{4} : S_{6}$ & $11520$ & $2$ & $4B,6A^2$\\
		$24$ & $2^{2+7}$ & $2^{2+6} : S_{3}$ & $1536$ & $2$ & $4B,4C^2$\\
		$25$ & $2^{5}$ & $2^{1+6} : \mbox{PSL}_{3}(2)$ & $21504$ & $2$ & $4C,6F^2$\\
		$26$ & $2^{2+4}$ & $2^{4} :3: S_{5}$ & $5760$ & $2$ & $4B,6C^2$\\
		$27$ & $2^{2}$ & $A_{8}$ & $20160$ & $2$ & $6A,10A^2$\\
		$28$ & $2^{6}$ & $2^{4} : A_{5}$ & $960$ & $2$ & $4B,10B^2$\\
		$29$ & $2^{1+6}$ & $2^{2+4} : S_{3}$ & $384$ & $2$ & $4B,6C^2$\\
		$30$ & $2^{5}$ & $2^{2+6} : S_{3}$ & $1536$ & $2$ & $4C,8B^2$\\
		$31$ & $2^{5}$ & $2^{2+6} : S_{3}$ & $1536$ & $2$ & $4C,10B^2$\\
		$32$ & $2^{1+5}$ & $S_{6}$ & $720$ & $2$ & $4B,10A^2$\\
		$33$ & $2^{6}$ & $S_{6}$ & $720$ & $2$ & $4B,12C^2$\\
		$34$ & $2^{4}$ & $2^{4} : \mbox{PSL}_{3}(2)$ & $2688$ & $2$ & $4C,6F^2$\\
		$35$ & $2^{1+3}$ & $2^{4} : \mbox{PSL}_{3}(2)$ & $2688$ & $2$ & $4C,10A^2$\\
		$36$ & $2$ & $A_{8}$ & $20160$ & $2$ & $6A,12C^2$\\
		$37$ & $2^{4}$ & $2^{6} : D_{10}$ & $640$ & $S_3$ & $8B^3$\\
		$38$ & $2^{6}$ & $2^{6} : S_{3}$ & $384$ & $2$ & $4B,8B^2$\\
		$39$ & $2^{1+4}$ & $2^{2+5} : S_{3}$ & $768$ & $2$ & $4C,6C^2$\\
		$40$ & $2$ & $2^{4} : S_{6}$ & $11520$ & $2$ & $6A,12C^2$\\
		$41$ & $2^{6}$ & $2\times \mbox{PSL}_{3}(2)$ & $336$ & $2$ & $4B,12C^2$\\
		$42$ & $2^{4}$ & $2^{6} : (S_{3} \times S_{3})$ & $2304$ & $1$ & $6A,6C,8B$\\
		$43$ & $2$ & $2^{4} : A_{6}$ & $5760$ & $2$ & $6A,8B^2$\\
		$44$ & $2^{1+0}$ & $2^{4} : S_{5}$ & $1920$ & $S_3$ & $6C^3$\\
		$45$ & $2^{2}$ & $S_{6}$ & $720$ & $S_3$ & $6C^3$\\
		$46$ & $1$ & $M_{11}$ & $7920$ & $2$ & $8B^2,10A$\\
		$47$ & $2^{1+0}$ & $2^{5} : S_{5}$ & $3840$ & $2$ & $6C,10A^2$\\
		$48$ & $1$ & $A_{7}$ & $2520$ & $S_3$ & $10A^3$\\
		$49$ & $2^{4}$ & $2^{3+4} : S_{3}$ & $768$ & $1$ & $6C,6C,8B$\\
		$50$ & $2^{1+2}$ & $S_{6}$ & $720$ & $2$ & $6C^2,10A$\\
		$51$ & $1$ & $2^{4} : \mbox{PSL}_{3}(2)$ & $2688$ & $3$ & $12C^3$\\
		$52$ & $2$ & $2^{4} : S_{5}$ & $1920$ & $2$ & $6C,8B^2$\\
		$53$ & $2^{4}$ & $2^{2+4} :3$ & $192$ & $2$ & $4C,10B^2$\\
		$54$ & $2^{4}$ & $2^{5} : S_{3}$ & $192$ & $2$ & $4C,12C^2$\\
		$55$ & $2^{2}$ & $2^{4} : (S_{3} \times S_{3})$ & $576$ & $2$ & $6C,10A^2$\\
		$56$ & $2$ & $2^{4} : S_{5}$ & $1920$ & $1$ & $8B,10A,12C$\\
		$57$ & $2$ & $2^{2+5} : S_{3}$ & $768$ & $2$ & $6C,12C^2$\\
		$58$ & $2$ & $2^{2+6} :3$ & $768$ & $2$ & $6F,10B^2$\\
		$59$ & $2$ & $2^{4} : (S_{3} \times S_{3}) :2$ & $1152$ & $1$ & $8B,10A,12C$\\
		$60$ & $2^{2}$ & $2^{2+3} : S_{3}$ & $192$ & $2$ & $8B,10A^2$\\
		$61$ & $2^{3}$ & $S_{3} \times S_{4}$ & $144$ & $1$ & $6C,10A,12C$\\
		$62$ & $1$ & $2^{2+3} : S_{3}$ & $192$ & $S_3$ & $10B^3$\\
		$63$ & $2$ & $2^{1+2+3} : S_{3}$ & $384$ & $1$ & $8B,10A,10B$\\
		$64$ & $2$ & $2^{1+4} : S_{3}$ & $192$ & $2$ & $6C,12C^2$\\
		$65$ & $2$ & $2^{5} : S_{3}$ & $192$ & $2$ & $8B,12C^2$\\
		$66$ & $1$ & $\mbox{GL}_{2}(4) :2$ & $360$ & $2$ & $10A,12C^2$\\
		$67$ & $1$ & $2^{5} : D_{10}$ & $320$ & $2$ & $10B^2,12C$\\
		$68$ & $2^{2}$ & $2^{4} : S_{3}$ & $96$ & $1$ & $6C,8B,12C$\\
		$69$ & $2^{2}$ & $2^{4} : S_{3}$ & $96$ & $1$ & $6C,10B,12C$\\
		$70$ & $1$ & $2^{2+4} : S_{3}$ & $384$ & $1$ & $8B,8B,12C$\\
		$71$ & $2$ & $2^{4} : S_{3}$ & $96$ & $2$ & $6C,8B^2$\\
		$72$ & $2$ & $2^{4} : S_{3}$ & $96$ & $2$ & $6C,10B^2$\\
		$73$ & $1$ & $2^{2+4} :3$ & $192$ & $2$ & $6F^2,8B$\\
		$74$ & $1$ & $2^{2+3} : S_{3}$ & $192$ & $2$ & $6F,12C^2$\\
		$75$ & $1$ & $2^{2+3} : S_{3}$ & $192$ & $2$ & $10B,12C^2$\\
		$76$ & $1$ & $\mbox{PSL}_{3}(2)$ & $168$ & $2$ & $10A,12C^2$\\
		$77$ & $2$ & $2^{2+4}$ & $64$ & $2$ & $10A,12C^2$\\
		$78$ & $2$ & $2^{2+4}$ & $64$ & $2$ & $10B,12C^2$\\
		$79$ & $2$ & $A_{5}$ & $60$ & $2$ & $6C,10B^2$\\
		$80$ & $2$ & $A_{5}$ & $60$ & $2$ & $10A,10B^2$\\
		$81$ & $1$ & $S_{5}$ & $120$ & $2$ & $8B^2,10A$\\
		$82$ & $1$ & $S_{5}$ & $120$ & $2$ & $10A,12C^2$\\
		$83$ & $2$ & $2^{4} : S_{3}$ & $96$ & $1$ & $8B,10A,12C$\\
		$84$ & $2$ & $2^{4} : S_{3}$ & $96$ & $1$ & $10A,10B,12C$\\
		$85$ & $1$ & $2^{5} : S_{3}$ & $192$ & $1$ & $12C,12C,12C$\\
		$86$ & $1$ & $2^{5} :3$ & $96$ & $2$ & $8B,12C^2$\\
		$87$ & $1$ & $2^{5}$ & $32$ & $S_3$ & $8B^3$\\
		$88$ & $1$ & $2^{1+4}$ & $32$ & $S_3$ & $10B^3$\\
		$89$ & $1$ & $(S_{3} \times S_{3}) :2$ & $72$ & $2$ & $8B,12C^2$\\
		$90$ & $1$ & $2^{2+2+3}$ & $128$ & $1$ & $8B,10B,12C$\\
		$91$ & $1$ & $2^{2+2+3}$ & $128$ & $1$ & $8B,10B,12C$\\
		$92$ & $2$ & $2^{2+3}$ & $32$ & $2$ & $6C,12C^2$\\
		$93$ & $1$ & $2^{3+3}$ & $64$ & $2$ & $12C,12C^2$\\
		$94$ & $1$ & $A_{5}$ & $60$ & $2$ & $8B^2,12C$\\
		$95$ & $1$ & $2^{1+2} : S_{3}$ & $48$ & $2$ & $10A,10B^2$\\
		$96$ & $1$ & $S_{3} \times S_{3}$ & $36$ & $2$ & $8B^2,12C$\\
		$97$ & $1$ & $A_{4}$ & $12$ & $S_3$ & $12C^3$\\
		$98$ & $1$ & $2^{3+2}$ & $32$ & $2$ & $6F,10B^2$\\
		$99$ & $1$ & $2^{3+2}$ & $32$ & $2$ & $8B,10B^2$\\
		$100$ & $1$ & $2^{5}$ & $32$ & $2$ & $8B,12C^2$\\
		$101$ & $1$ & $2^{1+4}$ & $32$ & $2$ & $10B,12C^2$\\
		$102$ & $1$ & $2^{1+2} : S_{3}$ & $48$ & $1$ & $8B,8B,10B$\\
		$103$ & $1$ & $2^{3}$ & $8$ & $S_3$ & $12C^3$\\
		$104$ & $1$ & $7:3$ & $21$ & $2$ & $6F,12C^2$\\
		$105$ & $1$ & $2\times A_{4}$ & $24$ & $1$ & $6F,10B,12C$\\
		$106$ & $1$ & $S_{4}$ & $24$ & $1$ & $8B,12C,12C$\\
		$107$ & $1$ & $S_{4}$ & $24$ & $1$ & $8B,12C,12C$\\
		$108$ & $1$ & $D_{12}$ & $12$ & $2$ & $12C,12C^2$\\
		$109$ & $1$ & $2^{2}$ & $4$ & $S_3$ & $12C^3$\\
		$110$ & $1$ & $D_{10}$ & $10$ & $2$ & $8B,10B^2$\\
		$111$ & $1$ & $3$ & $3$ & $S_3$ & $10B^3$\\
		$112$ & $1$ & $2^{3}$ & $8$ & $2$ & $10B,10B^2$\\
		$113$ & $1$ & $2^{3}$ & $8$ & $2$ & $10B^2,12C$\\
		$114$ & $1$ & $D_{12}$ & $12$ & $1$ & $8B,10B,12C$\\
		$115$ & $1$ & $S_{3}$ & $6$ & $2$ & $10B^2,12C$\\
		$116$ & $1$ & $S_{3}$ & $6$ & $2$ & $12C,12C^2$\\
		$117$ & $1$ & $D_{10}$ & $10$ & $1$ & $8B,10B,12C$\\
		$118$ & $1$ & $2^{3}$ & $8$ & $1$ & $10B,10B,12C$\\
		$119$ & $1$ & $D_{8}$ & $8$ & $1$ & $10B,12C,12C$\\
		$120$ & $1$ & $2^{3}$ & $8$ & $1$ & $10B,12C,12C$\\
		$121$ & $1$ & $2^{2}$ & $4$ & $2$ & $10B^2,12C$\\
		$122$ & $1$ & $S_{3}$ & $6$ & $1$ & $10B,12C,12C$\\
		$123$ & $1$ & $3$ & $3$ & $2$ & $10B,12C^2$
	\end{longtable}
}

%\addtocounter{table}{-1}  %  This is necessary, but don't know why 

%%%%%%%%%%%%%%%%%%%%%%%%%%%%%%%%%%%%%%

\subsection{Orbits on the feasible axes in $\Griess$ under $2.B \cap \N0$}
\label{N0capH:orbits}

Recall that the group $\M_{v^+}$ (of structure $2.B$)
is the stabilizer of the specific axis $v^+$, and that an axis is
{\em feasible} if it is in the orbit of the specific axis $v^-$ under $\Mv$.
In Table~\ref{tab:N0:2B:orbits} we list the $32$ orbits
on the feasible axes under the action of the
group $\M_{v^+} \cap \N0\cong 2^{2+11+20}.(M_{22}:2).S_3$.

\medskip

For each $(\M_{v^+} \cap \N0)$-orbit on the feasible axes, we display
the structure of the stabilizer in $\M_{v^+} \cap \N0$ of an axis in 
that orbit.
This stabilizer is presented in the form $C = E.G.S$, where
$E = C \cap 2^{2+11+20}$, $E.G = C \cap \M_{v^+} \cap \Nxyz$.
The groups 
$E$, $G$, and $S$ are displayed in Columns 2--5 of the table, 
with organization and legend as in Table~\ref{tab:N0:orbits}.

We have $\M_{v^+} \cap \Nxyz \vartriangleleft \M_{v^+} \cap \N0$ and
$(\M_{v^+} \cap \N0) / (\M_{v^+} \cap \Nxyz) = \N0 / \Nxyz = S_3$.
So the discussion of the decomposition of the $(\M_{v^+} \cap \N0)$-orbits
into $(\M_{v^+} \cap \Nx0)$-orbits and $(\M_{v^+} \cap \Nxyz)$-orbits
can be taken verbatim from the corresponding discussion in
Subsection~\ref{N0:orbits}. We list the fusion of
$(\M_{v^+} \cap \Nx0)$-orbits into $(\M_{v^+} \cap \C)$-orbits in
Column~6 of Table~\ref{tab:N0:2B:orbits} in the same format as in the
corresponding column in Table~\ref{tab:N0:orbits}.  
Here the names of the $(\M_{v^+} \cap \C)$-orbits are taken from
Table~\ref{table:orbits:H_axes}.

\smallskip

\ifdim \the\dimexpr\pagegoal-\pagetotal\relax < 20ex
\pagebreak
\fi

\begin{center}
	\captionof{table}{Orbits on feasible axes under $\M_{v^+} \cap {\N0} $}
	\label{tab:N0:2B:orbits}
\end{center}
\nopagebreak
\vspace{-5ex}
\nopagebreak
{\small 
	\begin{longtable}{|r|l|l|r|c|l|}
		\hline
		\multicolumn{6}{|c|}
		{$(\M_{v^+} \cap N_0)$-orbits on feasible axes, and stabilizers
			$\vphantom{2^{2^2}}C = E.G.S$  of these axes} \\
		\hline
		\textbf{No} & $E$ & $G$ & $|G|$ & $S$ & $(\M_{v^+} \cap \C)$-orbits \\
		\hline
		\endfirsthead
		
		\hline
		\textbf{No} & $E$ & $G$ & $|G|$ & $S$ & $(\M_{v^+} \cap \C)$-orbits \\
		\hline
		\endhead
		
		\hline
		\endfoot
		
		\hline
		\endlastfoot
		$\vphantom{2^{2^2}}$
$1$ & $2^{2+11+20} \hphantom{i}$ & $M_{22} :2$ & \hphantom{i} $887040$ & $2$ & $2A1,2A0^2$\\
$2$ & $2^{2+11+18}$ & $2^{5} : S_{5}$ & $3840$ & $S_3$ & $2A0^3$\\
$3$ & $2^{2+10+15}$ & $2^{4} : S_{6}$ & $11520$ & $2$ & $2A0,2B1^2$\\
$4$ & $2^{2+10+14}$ & $2^{4} : \mbox{PSL}_{3}(2)$ & $2688$ & $2$ & $2A0,2B0^2$\\
$5$ & $2^{1+11+10}$ & $\mbox{PSL}_{3}(4) :2$ & $40320$ & $2$ & $2A0,4A1^2$\\
$6$ & $2^{2+7+14}$ & $2^{4} : \mbox{PSL}_{3}(2)$ & $2688$ & $2$ & $2B1,2B0^2$\\
$7$ & $2^{2+7+12}$ & $2^{1+3+3} : S_{3}$ & $768$ & $S_3$ & $2B0^3$\\
$8$ & $2^{1+7+12}$ & $2^{5} : S_{5}$ & $3840$ & $2$ & $2B1,4A1^2$\\
$9$ & $2^{1+5+11}$ & $2^{4} : S_{6}$ & $11520$ & $2$ & $2B1,4B1^2$\\
$10$ & $2^{1+7+11}$ & $2^{2+5} : S_{3}$ & $768$ & $2$ & $2B0,4A1^2$\\
$11$ & $2^{6+10}$ & $2^{4} : S_{5}$ & $1920$ & $S_3$ & $4A1^3$\\
$12$ & $2^{1+4+10}$ & $2^{4} : \mbox{PSL}_{3}(2)$ & $2688$ & $2$ & $2B0,4C1^2$\\
$13$ & $2^{2+9}$ & $\mbox{PSL}_{3}(4) :2$ & $40320$ & $2$ & $4A1,6A1^2$\\
$14$ & $2^{1+5+10}$ & $2^{2+5} : S_{3}$ & $768$ & $2$ & $2B0,4B1^2$\\
$15$ & $2^{6+6}$ & $2^{4} : S_{6}$ & $11520$ & $2$ & $4A1,4B1^2$\\
$16$ & $2^{2+8}$ & $2^{5} : S_{5}$ & $3840$ & $S_3$ & $4B1^3$\\
$17$ & $2^{4+8}$ & $2^{4} : \mbox{PSL}_{3}(2)$ & $2688$ & $2$ & $4A1,4C1^2$\\
$18$ & $2^{4+10}$ & $2^{3+3} : S_{3}$ & $384$ & $2$ & $4A1,4B1^2$\\
$19$ & $2^{2+9}$ & $2^{4} : S_{5}$ & $1920$ & $2$ & $4A1,6C1^2$\\
$20$ & $2^{1+0}$ & $M_{22} :2$ & $887040$ & $2$ & $6A1,10A1^2$\\
$21$ & $2^{1+7}$ & $2^{4} : S_{6}$ & $11520$ & $1$ & $4B1,6A1,6C1$\\
$22$ & $2^{3+8}$ & $2^{2+4} : S_{3}$ & $384$ & $2$ & $4B1^2,4C1$\\
$23$ & $2^{3+6}$ & $2^{3+3} : S_{3}$ & $384$ & $S_3$ & $4B1^3$\\
$24$ & $2^{2+7}$ & $2^{2+5} : S_{3}$ & $768$ & $2$ & $4B1,4C1^2$\\
$25$ & $2^{2+4}$ & $2^{4} : S_{5}$ & $1920$ & $2$ & $4B1,6C1^2$\\
$26$ & $2^{1+6}$ & $2^{2+4} : S_{3}$ & $384$ & $2$ & $4B1,6C1^2$\\
$27$ & $2^{1+5}$ & $S_{6}$ & $720$ & $2$ & $4B1,10A1^2$\\
$28$ & $2^{1+3}$ & $2^{4} : \mbox{PSL}_{3}(2)$ & $2688$ & $2$ & $4C1,10A1^2$\\
$29$ & $2^{1+4}$ & $2^{2+5} : S_{3}$ & $768$ & $2$ & $4C1,6C1^2$\\
$30$ & $2^{1+0}$ & $2^{4} : S_{5}$ & $1920$ & $S_3$ & $6C1^3$\\
$31$ & $2^{1+0}$ & $2^{5} : S_{5}$ & $3840$ & $2$ & $6C1,10A1^2$\\
$32$ & $2^{1+2}$ & $S_{6}$ & $720$ & $2$ & $6C1^2,10A1$
	\end{longtable}
}

%%%%%%%%%%%%%%%%%%%%%%%%%%%%%%%%%%%%%%%%%%%%%%%%%%%%%%%%%%%%%%%%%%

\section{A certificate for computing the order of the Monster}
\label{appendix:certificate}

A key step in computing the order of the Monster is the enumeration
of the $\C$-orbits on the axes and their sizes, as listed
in Table~\ref{table:orbits:G_x0_axes}. These data can be readily
computed using a suitable linear algebra program, provided that
Table~\ref{table:action:G_x0_axes} has been established as correct.
However, generating the entries in Table~\ref{table:action:G_x0_axes}
involves sophisticated computations within the Monster group,
which are carried out using the mmgroup package.

Once these data (along with certain auxiliary information) are
available, their correctness can be verified through significantly
simpler computations. For this reason, it is practical to store all
the information required to generate Table~\ref{table:action:G_x0_axes}
in a separate human-readable file, enabling a simple and efficient
verification of its correctness. We refer to such a file as a
\emph{certificate}.

Similarly, we need the enumeration of the $(H \cap \Nxyz)$-orbits
on the feasible 2A axes and their sizes, as listed
in Table~\ref{table:orbits:H_axes}. These data can be readily
computed from Table~\ref{table:action:tau:H}. We therefore also 
compute a certificate for verifying the correctness of that
table.  

Our accompanying program code computes these certificates and stores
them in the subdirectory {\tt axis\_orbits/certificates}.
That subdirectory also contains Python scripts for verifying the
certificates. Therefore, it suffices to check the programs and data in this
subdirectory in order to verify the correctness of Tables~\ref{table:action:G_x0_axes} and~\ref{table:action:tau:H}.
For further details we refer the reader to the program code.
In the remainder of this section we focus on the verification of the
certificate for Table~\ref{table:action:G_x0_axes}. The verification of
the other certificate is similar. 

For verifying a certificate, we need only a tiny fraction of the
functionality of the mmgroup package. 
For each $\C$-orbit on the axes we store the following information
in a certificate:
\begin{enumerate}[itemsep=1pt, parsep=0pt]
	\item[-]
	The name of the orbit.
	\item[-]
	An element of $\M$ mapping the standard axis $v^+$ to a
	representative $v_i$ of the orbit.
	\item[-]
	A generating system for the stabilizer of axis $v_i$ in $\C$.
	\item[-]
	A list of elements of $\C$ mapping $v_i$ to a set of representatives
	of $\Nx0$-orbits contained in this orbit.
\end{enumerate}
Note that mmgroup allows storing elements of $\M$
in a compact form.
To compute the $\Nx0$-orbits (and their relative sizes) contained in
the $\C$-orbit of an axis $v_i$, we analyze the action of the stabilizer
of $v_i$ in $\C$ on the type-4 vectors in $\Lambda / 2\Lambda$, as
described in Section~\ref{subsection:speedup:axes}. For each case, we
verify that the resulting list of $\Nx0$-orbits agrees with the
corresponding list of precomputed $\Nx0$-orbit representatives in the
certificate.
For any such $\Nx0$-orbit representative $v_j$ we store the
following information in the certificate:
\begin{enumerate}[itemsep=1pt, parsep=0pt]
	\item[-]
	Orbit size (relative to $|X^+| : |\C : \Nx0|$).
	\item[-]
	Elements of $\C$ mapping the axes $v_j \cdot \tau^e$
	to the representatives of their $\C$-orbits, for $e=\pm 1$.
\end{enumerate}
Given such a certificate, we may compute 
Table~\ref{table:action:tau:H} as described in 
Sections~\ref{subsection:basic:axes}--\ref{subsection:speedup:axes},
but with the following simplifications:
\begin{enumerate}[itemsep=1pt, parsep=0pt]
	\item[-]
	If a representative of a $\C$- or $\Nx0$-orbit is required,
	we may compute it from the data in the certificate.
	\item[-]
	If an axis must be transformed to the representative of its
	$\C$-orbit, we may read the transformation from the certificate.
\end{enumerate}
For this computation, it is sufficient to be able to perform the following
operations related to the Monster with mmgroup: 
\begin{enumerate}[itemsep=1pt, parsep=0pt]
	\item[-]
	Transforming an axis in $\Griess$ by an element of $\M$.
	\item[-]
	Mapping an element of $\C$ to the automorphism group
	$\Co_1$ of $\Lambda / 2 \Lambda$.
	\item[-]
	Computing in the natural representation of $\Co_1$ on
	$\Lambda / 2 \Lambda$. 
\end{enumerate}
In particular, we do not need any word-shortening algorithm for the
Monster. During such a computation we also check the internal
consistency of the certificate. The size of the certificate is about
60~kBytes; its verification takes about 20~seconds on the authors' PCs.

%%%%%%%%%%%%%%%%%%%%%%%%%%%%%%%%%%%%%%%%%%%%%

\section{Classes of involutions in the Monster}
\label{appendix:classes:involutions:Monster}

It is well known that the Monster has two classes of involutions
called 2A and 2B in~\cite{Atlas}. 
In this appendix we will prove this
fact using the character table of the group $\mbox{Co}_1$ and
computations in mmgroup only. 

\medskip

From Theorem~\ref{theorem:order:monster} it follows that $\C$ has odd index in $\M$.
Thus it suffices to show that all involutions in $\C$ fuse either to class 2A or
to class 2B in the Monster.

The Python script {\tt involutions\_G\_x0.py} in the accompanying
program code performs a brute-force calculation with mmgroup to
establish this fact. It turns out that $\C$ has seven classes
of involutions, two classes fusing to class 2A in $\M$, and
five classes fusing to class 2B in $\M$.

\medskip

Alternatively, this last fact can be established using the character information about $\M$ and $\C$
stored in the GAP computer algebra package~\cite{GAP4} with the following GAP program:
 
\vspace{0.5cm}
 
\hrule
{\small
\begin{verbatim}
LoadPackage("ctbllib");;
chM := CharacterTable("M");;
chGx0 := CharacterTable("2^1+24.Co1");;
fusion := FusionConjugacyClasses(chGx0, chM);;
involutions := Positions(OrdersClassRepresentatives(chGx0), 2);;
List(involutions, t->[ClassNames(chGx0)[t], ClassNames(chM)[fusion[t]]]);
\end{verbatim}
}
\hrule

\end{appendices}

%%%%%%%%%%%%%%%%%%%%%%%%%%%%%%%%%%%%%%%%%%%%%%%%%%%%%%%%%%%%%%%%%%%%%%%%%%%%%%

%\bibliography{references}{}
%\bibliographystyle{alpha}
\newcommand{\etalchar}[1]{$^{#1}$}

\end{document}